\documentclass[a4paper]{article}
\usepackage{amssymb}
\usepackage[german,english]{babel}
\usepackage{theorem}
\usepackage{QuaAnDefs}

\MyPage{16cm}{23.5cm}{0cm}{FRpage} 
\MyHead{{\sl Kaup, Reiffen: Quasi-analytische Zerlegungen}}

\SecNumStyle{2}
\MyFoot{}
\def\4#1{{\bf #1}}
\def\TT{{\rm T}}

\def\Rang{\mathop{\rm Rang}\nolimits}
\def\Reg{\mathop{\rm Reg}\nolimits}

\def\dist{\mathop{\rm dist}\nolimits}
\Itemiarabic

\def\Bul#1{#1_{\bullet}}  
\def\0#1{#1_{0}}  
\def\1#1{#1_{1}}  
\def\2#1{#1_{2}}  
\def\IND#1#2{#1_{#2}} 
\def\dimMin{\mathop{\underline{\rm dim}}\nolimits}
\def\calA{{\cal A}}

\def\rank{\mathop{\rm rank}\nolimits}
\def\DEFI#1{\DEF{#1}\index{#1}}
\def\AA{{\bf A}}
\def\DF{\5D_{\5F}}
\def\Df{\5D_{f}}
\let \Dot = \Bul

\begin{document}

\centerline{\bf\huge  Quasi-analytische Zerlegungen}

\bigskip

\centerline{\sc Burchard Kaup und Hans-J{\"o}rg Reiffen}

\null
\bigskip

\selectlanguage{english}
\begin{abstract}
   The leaves in singular holomorphic foliation theory are examples of quasi-analytic layers. In the first part of our publication we are concerned with a theory of these subjects. A quasi-analytic decomposition of a complex manifold is a decomposition into pairwise disjoint connected quasi-analytic layers. These are holomorphic foliations in the sense of P. Stefan and K. Spallek.  A very different but more usual conception of holomorphic foliations is develloped by P.~Baum and R.~Bott. It is based on holomorphic sheaf theory. In the second part we study the relation between quasi-analytic decompositions and singular holomorphic foliations in the sense of Baum and Bott.
\end{abstract}

\bigskip

\centerline{2000 Mathematics Subject Classification: 32C15, 32S65.}

\selectlanguage{german}

\subsection{Einleitung}\label{Ein}

Es sei $X$ eine zusammenh\"angende komplexe Mannigfaltigkeit mit abz\"ahlbarer Topologie und $Z \subset X$ eine Teilmenge von $X$. In der Theorie der holomorphen Bl\"atterungen sind solche Teilmengen $Z$ von Interesse, zu denen es eine injektive holomorphe Immersion $\Mapping \phi Y X$  eines komplexen Raumes $Y$ in $X$ mit $Z$ als Bild gibt. Dann induziert $\phi$ auf $Z$ eine komplexe Struktur. Wir nennen eine solche Struktur \DEF{quasi-analytisch} und $Z$ mit dieser Struktur eine \DEF{Schichtung}. Ist $Z$ bez\"uglich der Struktur zusammenh\"angend, so nennen wir $Z$ eine \DEF{Schicht}. 
Die Bl\"atter in einer holomorphen Bl\"atterungstheorie sind stets Schichten. 

In Abschnitt \ref{QaS}
geben wir eine Definition der quasi-analytischen Struktur, ohne holomorphe Abbildungen zu benutzen. Eine quasi-analytische Struktur auf $Z$ ist eine Topologie $\5T$ auf $Z$ derart,  da{\ss} $(Z,\5T)$ lokal 
bez\"uglich $\5T$ jeweils eine lokal-analytische Teilmenge von $X$ ist (vgl.~Definition \ref{1.1}). Dadurch ist in nat\"urlicher Weise auch eine komplexe Struktur auf $(Z,\5T)$ festgelegt. Wir stellen in Abschnitt \ref{QaS} einige grundlegende Begriffe und Aussagen zusammen. Wie Beispiele in Abschnitt \ref{Bsp} zeigen, sind quasi-analytische Strukturen, auch wenn es sich um Schichten handelt, nicht eindeutig bestimmt. In Abschnitt \ref{AI}
machen wir unter der Einschr\"ankung einer abz\"ahlbaren Topologie (Schichten erf\"ullen diese Voraussetzung) eine abgeschw\"achte Eindeutigkeitsaussage (vgl.~\ref{3.10} und \ref{3.11}). 
Wie wir in Abschnitt \ref{QaS} zeigen, gilt allerdings eine scharfe Eindeutigkeitsaussage f\"ur Wege-vertr\"agliche quasi-analytische Strukturen (vgl.~\ref{1.6}). Dabei hei{\ss}t eine quasi-analytische Struktur \DEF{Wege-vertr\"aglich}, wenn jeder $X$-Weg bereits ein $\5T$-Weg ist. Die Existenz einer solchen Struktur auf $Z$ h\"angt allerdings von der Relativtopologie von $Z$ bzgl.~$X$ ab (vgl. \ref{1.5}.\ref{1.5.3}). 
Es handelt sich, wenn es sie gibt, um eine nat\"urliche quasi-analytische Struktur. In Abschnitt \ref{AI} gehen wir allgemein auf nat\"urliche quasi-analytische Strukturen ein; wir nennen sie \DEF{schwach-analytische Strukturen}.
Man konstruiert sie mit Hilfe des sogenannten analytischen Inhalts. Das ist die Menge aller Teilmengen $A$ von $Z$, die lokal-analytische Teilmengen von $X$ sind. Nicht jede Teilmenge $Z$ tr\"agt eine schwach-analytische Struktur. Wenn eine solche existiert, ist sie eindeutig bestimmt (vgl.~\ref{3.19}).

Seit Beginn der 1970er Jahre studiert man holomorphe Bl\"atterungen mit Singularit\"aten. 1972 hat {\sc H.~Holmann} einen speziellen Typ von singul\"aren holomorphen Bl\"atterungen mit Hilfe gewisser holomorpher Abbildungen eingef\"uhrt (vgl.~\cite{Holmann}). Er konnte dabei wie im Fall regul\"arer Bl\"atterungen einen Blattbegriff einf\"uhren. Problematisch ist dies bei dem ebenfalls in 1972 von {\sc P.~Baum} und {\sc R.~Bott} eingef\"uhrten Begriff einer singul\"aren holomorphen Bl\"atterung mit Hilfe von Garben holomorpher Vektorfelder bzw.~Pfaffscher Formen 
(vgl.~\cite{BaumBott} und \cite{Baum}).
Dieser Begriff besitzt jedoch den gew\"unschten Grad von Allgemeinheit. 

Deshalb gehen wir von ihm in dieser Arbeit aus. In \cite{Reiffen1} wurde f\"ur derartige Bl\"atterungen ein Blattbegriff eingef\"uhrt, welcher den von {\sc H.~Holmann} verallgemeinert. Allerdings besitzt nicht jede singul\"are holomorphe 
Bl\"atterung Bl\"atter \"uberall. Wir beziehen uns in dieser Arbeit bei Verweisen zur Bl\"atterungstheorie stets auf 
\cite{Reiffen2}, eine Arbeit, die \"uber das Internet erreichbar ist. 

Inspiriert durch Arbeiten von {\sc P.~Stefan} hat {\sc K.~Spallek} einen Bl\"atterungsbegriff f\"ur differenzierbare R\"aume eingef\"uhrt (vgl.~\cite{Spallek}). F\"ur den holomorphen Fall und komplexe Mannigfaltigkeiten als Tr\"agerraum stimmt sein Begriff mit dem von uns in  Abschnitt \ref{QaZ} definierten Begriff der quasi-analytischen Zerlegung \"uberein. Eine quasi-analytische Zerlegung von $X$ ist eine quasi-analytische Struktur auf $Z = X$; die zugeh\"origen Schichten nennen wir \DEF{Bl\"atter}. 

In den Abschnitten  \ref{AdZ} bis \ref{KZ} geht es um den Zusammenhang zwischen quasi-analytischen Zerlegungen und holomorphen Bl\"atterungen.

In Abschnitt \ref{QaZ} stellen wir allgemeine Grundbegriffe f\"ur die Theorie der quasi-analytischen Zerlegungen zusammen und gehen vor allem auf zwei zentrale Beispiele ein: Ist $\Theta'$ eine koh\"arente involutive Untergarbe der Garbe $\Theta$ der holomorphen Vektorfelder auf $X$, so definiert $\Theta'$, wie Spallek gezeigt hat, eine glatte quasi-analytische Zerlegung $\5D'$ von $X$. Wir nennen sie die \DEF{Spallek-Zerlegung} zu $\Theta'$; sie ist insbesondere Wege-vertr\"aglich (vgl. \ref{4.8} und \ref{4.9}). Ist $\5F$ eine holomorphe Bl\"atterung auf $X$ mit Bl\"attern \"uberall, so definieren diese eine quasi-analytische Zerlegung $\5D_{\5F}$ von $X$; sie ist schwach-analytisch (vgl.~\ref{4.13}). Wir behandeln die Beziehung zwischen $\5D_{\5F}$ und der gleichfalls existierenden Spallek-Zerlegung $\5D_{\5F}'$ von $\5F$ (vgl.~\ref{4.16}). In Abschnitt \ref{AdZ} behandeln wir quasi-analytische Zerlegungen und holomorphe Bl\"atterungen, die abbildungsdefiniert sind, d.h.~lokal durch holomorphe Abbildungen beschrieben werden k\"onnen (vgl. \ref{5.3}, \ref{5.5}). Der Begriff der Regularit\"at einer quasi-analytischen Zerlegung $\5D$ wird in naheliegender Weise eingef\"uhrt (vgl. \ref{5.6}) und damit der der Singularit\"atenmenge $\Sing \5D$. F\"ur unsere Theorie wichtig ist der Begriff der fasertreuen holomorphen Abbildung (vgl. \ref{5.9}); insbesondere sind offene Abbildungen fasertreu. Wir erhalten das Ergebnis: die abbildungsdefinierten quasi-analytischen Zerlegungen mit fasertreuen lokalen Beschreibungen entsprechen genau den abbildungsdefinierten holomorphen Bl\"atterungen mit fasertreuen lokalen Beschreibungen  (vgl. \ref{5.5}, \ref{5.15}). 
Wie Beispiele am Ende von Abschnitt \ref{AdZ} zeigen, k\"onnen quasi-analytische Zerlegungen sehr chaotisch sein, selbst wenn die zugeh\"orige Vektorfeldgarbe koh\"arent ist. Deshalb erscheint es sinnvoll, bei der Definition der 
Koh\"arenz einer quasi-analytischen Zerlegung $\5D$ nicht nur die Koh\"arenz der Vektorfeldgarbe $\Theta^{\5D}$ zu fordern, sondern dar\"uber hinaus eine st\"arkere Bindung von $\5D$ und $\Theta^{\5D}$ zu verlangen (vgl.~\ref{6.2}). In diesem Fall definiert $\5D$  in nat\"urlicher Weise eine holomorphe Bl\"atterung $\5F^{\5D}$ und auf 
$X \sm \Sing\5D$ stimmt die von dieser dort definierte Zerlegung mit $\5D|_{X \sm \Sing \5D}$ \"uberein. Die quasi-analytische Zerlegung $\5D$ hei{\ss}t vollst\"andig, wenn sie koh\"arent ist und wenn $\Theta^{\5D}$ vollst\"andig ist (vgl.~\ref{6.7}). Sie hei{\ss}t perfekt, wenn sie koh\"arent ist, wenn $\5F = \5F^{\5D}$ Bl\"atter \"uberall besitzt und wenn $\5D = \5D_{\5F}$ ist (vgl.~\ref{6.15}, \ref{6.16}). 

\"Ubrigens: unser Koh\"arenzbegriff ist funktionentheoretisch motiviert und stimmt nicht mit dem eher geometrisch motivierten Koh\"arenzbegriff von {\sc Spallek} \"uberein. 

Wir zeigen: 

\ref{6.18}: {\sl Die quasi-analytische Zerlegung $\5D$ sei abbildungsdefiniert mit fasertreuen lokalen Beschreibungen. Dann ist $\5D$ vollst\"andig und perfekt. }

\ref{6.19}: {\sl Die quasi-analytische Zerlegung $\5D$ sei koh\"arent, lokal-eigentlich, rein $p$-dimensional und es gelte 
$\dim \Sing \5D < p$. Dann ist $\5D$ vollst\"andig und perfekt.}

\ref{6.21}:
{\sl Die quasi-analytische Zerlegung $\5D$ sei lokal eigentlich und rein $1$-codimensional. Dann sind folgende Aussagen \"aquivalent:
\BegEN
   \item $\5D$ ist koh\"arent,
   \item $\5D$ ist vollst\"andig,
   \item $\5D$ ist perfekt,
   \item $\5D$ ist abbildungsdefiniert.\label{6.21.4a}
\EndEN
Im Fall von \ref{6.21.4a} sind die lokalen Beschreibungen \"ubrigens sofort schon fasertreu.}

In Abschnitt \ref{Anh}, dem Anhang, stellen wir einige Hilfs\"uberlegungen allgemeinerer Art zusammen.

Eine verk\"urzte Vorabversion dieser Arbeit wurde in \cite{Knoche} publiziert. Leider sind uns dort einige Fehler unterlaufen:
eine Schicht ist nicht immer Wege-vertr\"aglich, und
im Theorem 1.(1).(a) fehlt die Bedingung, da{\ss} $\5D$ lokal eigentlich ist (vgl.~Beispiel \ref{6.26} in diesem Text).

\subsection{Quasi-analytische Schichtungen und Schichten}\label{QaS}

Es sei $X$ stets eine $n$-dimensionale zusammenh\"angende komplexe Mannigfaltigkeit mit abz\"ahlbarer Topologie. Ist $M$ eine beliebige Teilmenge von $X$, dann bezeichnen wir die Relativtopologie von $M$ auch als die 
\DEFI{$X$-Topologie}.

\Beg{Definition}\label{1.1}
   Eine \DEFI{quasi-analytische Struktur} auf einer Teilmenge $Z$ von $X$ ist eine Topologie $\5T$ auf $Z$, f\"ur die gilt:
   \\
   Zu jedem $z \in Z$ gibt es eine offene $\5T$-Umgebung $A$ von $z$ in $Z$ und eine offene $X$-Umgebung 
   $U$ von $z$ in $X$ derart, da{\ss} $A$ eine analytische Teilmenge von $U$ ist und da{\ss} die $\5T$-Topologie  
   von $A$ mit der $X$-Topologie von $A$ \"ubereinstimmt.
   \\
   Wir nennen den topologischen Raum $(Z,\5T)$ eine \DEF{quasi-analytische Schichtung}\index{Schichtung}.
    Wenn in der obigen Bedingung $A$ und $U$ zusammenh\"angend sind (was man durch geeignete Verkleinerungen 
    stets erreichen kann), dann nennen wir $(A,U)$ ein \DEF{$\5T$-Pl\"attchen}\index{Pl\"attchen}.
   \\
   Die $\5T$-Zusammenhangskomponenten einer quasi-analytischen Schichtung (zusammen mit
   der durch $\5T$ induzierten Topologie) nennen wir 
   \DEF{quasi-analytische Schichten}.\index{Schicht}\index{quasi-analytische Schicht}
\End

Wir sprechen im Weiteren kurz von \DEFI{Schichtung} bzw.~\DEFI{Schicht}. Jede Schicht ist wieder eine Schichtung.

Da analytische Mengen lokal wegzusammenh\"angend sind, ist eine Schichtung $(Z,\5T)$ bez\"uglich $\5T$ ebenfalls lokal wegzusammenh\"angend und ihre Schichten sind genau ihre $\5T$-Wegzusammenhangskomponenten.

Die Topologie $\5T$ einer Schichtung ist feiner als die $X$-Topologie; lokal bzgl.~$\5T$ tr\"agt eine Schichtung die Relativtopologie von $X$. 

Ist $(Z,\5T)$ eine Schichtung in $X$, dann induziert $\5T$ in nat\"urlicher Weise eine komplexe Struktur auf $Z$:
ist $(A,U)$ ein $\5T$-Pl\"attchen, so ist $A$ mit seiner nat\"urlichen komplexen Struktur eine komplexe Karte
von $Z$. Wir denken uns eine Schichtung $(Z,\5T)$ stets mit dieser komplexen Struktur versehen. Es ist eine reduzierte komplexe Struktur.

Wir setzen im Weiteren alle betrachteten komplexen Strukturen als reduziert voraus.

Eine holomorphe Abbildung $\Mapping \phi Y {\Bul Y}$ zwischen komplexen R\"aumen $Y$ und $\Bul Y$ hei{\ss}t eine
\DEF{holomorphe Immersion}\index{Immersion}, wenn es zu jedem Punkt $y \in Y$ eine offene Umgebung $U$ von $y$ und eine offene Umgebung $\Bul U$ von $\Bul y = \phi(y)$ gibt derart, da{\ss} $\phi(U)$ eine analytische Teilmenge von $\Bul U$ und $\Mapping \phi U {\phi(U)}$ biholomorph ist.

\Beg{}\label{1.2}\
\BegEN
   \item \label{1.2.1}
      Sei $(Z,\5T)$ eine Schichtung in $X$. Dann ist die kanonische Inklusion $\Mapping \iota Z X$ eine injektive holomorphe 
      Immersion mit $\iota(Z) = Z$. 
   \item\label{1.2.2}
      Ist umgekehrt $\Mapping \phi Y X$ eine injektive holomorphe Immersion, so ist die Topologie $\5T$ auf $Z:=\phi(Y)$, 
      bez\"uglich
      der $\phi$ ein Hom\"oomorphismus ist, eine quasi-analytische Struktur auf $Z$, also $(Z,\5T)$ eine Schichtung. 
      $\Mapping \phi Y Z$ ist eine biholomorphe Abbildung.
\EndEN
\End

{\bf Beweis:} \ref{1.2.1} ist klar. --- \ref{1.2.2} ergibt sich sofort aus der Definition der holomorphen Immersion. \Qed

Eine holomorphe Abbildung wie in \ref{1.2}.\ref{1.2.2} hei{\ss}t eine \DEFI{Darstellung} der Schichtung $(Z,\5T)$. 
Wegen \ref{1.2}.\ref{1.2.1} besitzt jede Schichtung eine Darstellung. Deshalb sind Schichtungen genau die Teilmengen von $X$, die Bilder komplexer R\"aume bez\"uglich injektiver holomorpher Immersionen mit Zielraum 
$X$ sind.

\Beg{}\label{1.3}\
\BegEN
   \item\label{1.3.1}
      Sei $A$ eine lokal-analytische Teilmenge von $X$. Dann ist die $X$-Topologie von $A$ eine 
      quasi-analytische Struktur auf $A$,
      die sogenannte \DEFI{Standardstruktur}.
   \item\label{1.3.2}
      Zur Schichtung $(Z,\5T)$ in $X$ gebe es eine offene $X$-Umgebung $U$ von $Z$ derart, 
      da{\ss} f\"ur jede $X$-kompakte
      Teilmenge $K$ von $U$ die Menge $K \cap Z$ ebenfalls $\5T$-kompakt ist. 
      Dann ist $Z$ eine lokal-analytische Teilmenge von $X$ und
      $\5T$ ist die Standardstruktur, also die $X$-Topologie.
\EndEN
\End

{\bf Beweis} von \ref{1.3.2}: Sei $\Mapping \phi Y X$ eine Darstellung von $(Z,\5T)$. Dann ist $\Mapping \phi Y U$ eine eigentliche Abbildung. Aufgrund eines Satzes von Remmert (vgl.~\cite[45.17]{KK}) ist $Z = \phi(U)$ eine analytische Teilmenge von $U$, also eine lokal-analytische Teilmenge von $X$. Die Abbildung $\Mapping \phi Y Z$ (wobei $Z$ die $X$-Topologie trage) ist bijektiv, stetig und eigentlich, also ein Hom\"oomorphismus. Also ist $\5T$ die $X$-Topologie. \Qed

In \ref{1.3}.\ref{1.3.1} erf\"ullt $A$ mit der Standardstruktur die Voraussetzungen von \ref{1.3}.\ref{1.3.2}. Also gibt \ref{1.3} eine genaue Charakterisierung derjenigen Schichtungen in $X$ an, die lokal-analytische Teilmengen von $X$ mit ihrer Standardstruktur sind.

Die folgenden zwei Beispiele zeigen, wie vorsichtig man mit dem Begriff der Schichtung vorgehen sollte:

\Beg{}\label{1.4}\
\BegEN
   \item\label{1.4.1}
      Sei $M$ eine beliebige (nichtleere) Teilmenge von $\C$ (z.B.~$M = \7Q$, $\R \sm \7Q$, $\C$) und 
      $Y := \bigcup_{z \in M} \{z\} \times \C$. 
      Dann ist $Y$ in nat\"urlicher Weise eine eindimensionale komplexe Mannigfaltigkeit mit    
      eventuell \"uberabz\"ahlbar vielen Zusammenhangskomponenten. 
      Die nat\"urliche Inklusion $Y \hookrightarrow \C^2$ definiert eine quasi-analytische Struktur auf $Y \subset \C^2$.
   \item\label{1.4.2}
      F\"ur eine beliebige Teilmenge $Y \subset X$ ist $\tilde Y := \bigcup_{y \in Y} \{y\}$ mit der diskreten Topologie 
      auf nat\"urliche Weise
      eine null-dimensionale komplexe Mannigfaltigkeit. Die nat\"urliche Inklusion $\tilde Y \hookrightarrow X$ definiert eine 
      quasi-analytische Struktur auf $Y$.
\EndEN
\End

Weil $X$ nach Voraussetzung abz\"ahlbare Topologie hat und jede Schicht eine Darstellung besitzt, folgt auf Grund des Satzes von PoincarŽ-Volterra (vgl.~\cite{Bou}), da{\ss} jede Schicht eine abz\"ahlbare Topologie besitzt. Wir nennen eine Schichtung 
\DEF{abz\"ahlbar}\index{abz\"ahlbare Schichtung}\index{Schichtung!abz{\"a}hlbare}, 
wenn sie h\"ochstens abz\"ahlbare viele Schichten enth\"alt. Das ist genau dann der Fall, wenn die Schichtung abz\"ahlbare Topologie hat.

Eine lokal-analytische Teilmenge von $X$ mit der Standardstruktur ist eine abz\"ahlbare Schichtung. Im ¬\"Ubrigen verweisen wir auf die Beispiele \ref{1.4}.\ref{1.4.1}. 

\Beg{Satz und Definition}\label{1.5}
   Eine Schichtung $(Z,\5T)$ hei{\ss}t \DEFI{Wege-vertr\"aglich}, wenn die folgenden \"aquiva\-lenten Aussagen gelten:
   \BegEN
   \item\label{1.5.1}
      Jeder $X$-Weg $\Mapping \gamma {I = \Clint 01} Z$ ist ein $\5T$-Weg.
   \item\label{1.5.2}
      F\"ur jedes $\5T$-Pl\"attchen $(A,U)$ ist $A$ eine $X$-Wegzusammenhangskomponente von $Z \cap U$.
   \item \label{1.5.3}
      Zu jedem $z \in Z$ gibt es ein $\5T$-Pl\"attchen $(A,U)$ mit $z \in A$ derart, da{\ss} $A$ die 
      $X$-Wegzusammenhangs\-komponente
      von $Z \cap U$ mit $z \in A$ ist.
   \EndEN
\End

{\bf Beweis} \ref{1.5.1} $\IfThen$ \ref{1.5.2}: Da $A$ wegzusammenh\"angend ist, gen\"ugt es zu zeigen:
f\"ur jeden $X$-Weg
\\
$\Mapping \gamma {I = [0, 1]}{Z \cap U}$ mit  $\gamma(0) \in  A$  ist $J := \{t \in I :  \gamma(t) \in A\} = I$. 
Das ist aber klar: da $A$ abgeschlossen ist in $U$, ist $J$ abgeschlossen; da $\gamma$ 
nach Voraussetzung auch ein $\5T$-Weg ist und $A$ $\5T$-offen in $Z \cap U$, ist $J$ offen in $I$.
\\
\ref{1.5.3} $\IfThen$ \ref{1.5.1}: Sei $\gamma$ wie in \ref{1.5.1}; wir zeigen, da{\ss} $\gamma$  in allen $t_0 \in I$  
$\5T$-stetig   ist. Dazu sei $(A,U)$
ein $\5T$-Pl\"attchen mit $\gamma(t_0) \in A$ wie in \ref{1.5.3}. Weil $A$ eine $X$-Zusammenhangskomponente von 
$Z \cap U$ ist, ist 
$\gamma(t) \in A$ f\"ur $t$ nahe bei $t_0$. Weil $\5T$ die Relativtopologie von $A$ in $U$ ist, 
ist $\gamma$ auch
$\5T$-stetig in einer Umgebung von $t_0$.  \Qed

Wir sagen im Fall von \ref{1.5} auch, da{\ss} $\5T$ \DEFI{Wege-vertr\"aglich} ist.

\Beg{Beispiel}\label{Beisp1}
Ist $M$ in \ref{1.4}.\ref{1.4.1} total unzusammenh\"angend (z.B.~das Cantorsche Diskontinuum), dann ist die dadurch definierte Schichtung Wege-vertr{\"a}glich. F{\"u}r $M = \C$ ist das zugeh{\"o}rige
$Y$ nicht Wege-vertr{\"a}glich. 
\End

Weitere Beispiele zum Begriffe "Wege-vertr{\"a}glichÓ findet man im Abschnitt \ref{Bsp}.

\Beg{Satz}\label{1.6}
   Seien $\5T$ und $\Bul {\5T}$ zwei Wege-vertr\"agliche quasi-analytische Strukturen auf $Z$. Dann ist $\5T = \Bul{\5T}$.
\End

{\bf Beweis:}
Sei $z \in Z$. Dann gibt es zu $z$ ein $\5T$-Pl\"attchen $(A,U)$ und ein $\Bul{\5T}$-Pl\"attchen $(\Bul A,\Bul U)$ gem\"a{\ss} 
\ref{1.5.3}. Wir betrachten die $X$-Wegzusammenhangskomponente $\tilde A$ von $Z \cap U \cap \Bul U$, die $z$ enth\"alt. Weil
analytische Mengen lokal wegzusammenh\"angend sind, ist $\tilde A$ eine offene Teilmenge von $A$ und von $\Bul A$. Sei 
$\tilde U$ die $X$-Zusammenhangskomponente von $U \cap \Bul U$, die $z$ enth\"alt, dann ist $(\tilde A, \tilde U)$ ein
$\5T$- und ein $\Bul{\5T}$-Pl\"attchen. \Qed

\Beg{}\label{1.7}
   Sei $(Z,\5T)$ eine Schichtung und $V  \subset Z$ eine $\5T$-offene Teilmenge. 
   Dann ist $(V, \5T|_{V})$ wieder eine Schichtung, 
   die wir die \DEF{Einschr\"ankung}\index{Einschr\"ankung einer Schichtung} der Schichtung $(Z,\5T)$ nennen. 
\End

Es sei $(Z,\5T)$ eine Schichtung und $W \Opss X$, also $W$ eine offene Teilmenge von $X$. Dann ist $V := Z \cap W$ eine $\5T$-offene Teilmenge von $Z$. Die Schichten von $(V,\5T|_{V})$ nennen
wir kurz die \DEF{$\5T$-Komponenten}\index{T-Komponenten@$\5T$-Komponenten} von $Z \cap W$.

\Beg{Definition}\label{1.8}
   Die Schichtung $(Z,\5T)$ in $X$ hei{\ss}t \DEFI{ordentlich}\index{Schichtung!ordentliche}, 
   wenn es zu jedem Punkt $z \in Z$ eine $X$-Umgebungsbasis von
   $X$-offenen Umgebungen $W$ von $z$ gibt mit folgender Eigenschaft:
   \begin{quote}\parskip=0pt
   Seien $A,B$ zwei verschiedene $\5T$-Komponenten von $Z \cap W$, so gilt $\8A \cap B = \emptyset$, wobei $\8A$ der
   $W$-Abschlu{\ss} von $A$ sei. 
   \end{quote}
   {\parskip 0pt Sind f{\"u}r alle $W$ alle $\5T$-Komponenten von $W$ jeweils $W$-abgeschlossen, 
   so hei{\ss}t die Schichtung $(Z,\5T)$ \DEFI{sehr ordentlich}\index{Schichtung!sehr ordentliche}.}
\End

In der Situation \ref{1.8} nennen wir auch $\5T$ (sehr) ordentlich.

\Beg{Satz}\label{1.9}
   Die Schichtung $(Z,\5T)$ in $X$ sei abz\"ahlbar und ordentlich. Dann ist $(Z,\5T)$ Wege-vertr\"aglich.
\End

 {\bf Beweis:} Sei $z \in Z$ und $(A,U)$ ein  $\5T$-Pl\"attchen  mit $z \in A$. Sei $W \subset U$ eine offene
 $X$-Umgebung aus der in \ref{1.8} genannten Umgebungsbasis. Wir betrachten die $\5T$ Zusammenhangskomponente $A_{0}$ von
 $A \cap W$, die $z$ enth\"alt. $Z \cap W$ besitzt h\"ochstens abz\"ahlbar viele $\5T$-Komponenten $B_{\nu}$, $\nu \in N \subset \N$. Es ist $\8B_{\nu} \cap B_{\mu} = \emptyset$ f\"ur $\nu \not= \mu$. Sei $0 \in N$ und $z \in B_{0}$. Dann ist $A_{0}=B_{0}$. 
 Sei $\Mapping \gamma {\Clint 01}{Z \cap W}$ ein $X$-Weg mit $\gamma(0)=z$.
 F\"ur jedes $\nu \in N$ ist $\gamma^{-1}(B_{\nu}) = \gamma^{-1}(\8B_{\nu})$ eine abgeschlossene Teilmenge von $I$. Wegen Satz
 \ref{7.2} ist $\gamma^{-1}(B_{0})=I$, also $\gamma(I) \subset A_{0}$. Deshalb ist $A_{0}$ die $X$-Wegzusammenhangskomponente
 von $Z \cap W$, die $z$ enth\"alt. Nun kann mit $A_{0}$ die Situtation von $\ref{1.5}.\ref{1.5.3}$ hergestellt werden. \Qed
 
Ist im Beispiel \ref{1.4}.\ref{1.4.1} die Teilmenge $M \subset \C$ abz\"ahlbar, dann ist $Y$ eine abz\"ahlbare sehr ordentliche Schichtung, insbesondere also, wie in \ref{Beisp1} schon festgestellt, Wege-vertr\"aglich.

 Die Standardstruktur einer lokal-analytischen Teilmenge $A$ von $X$ ist sehr ordentlich, insbesondere Wege-vertr\"aglich.

 Der folgende Satz gibt ein einfaches Konstruktionsverfahren f\"ur quasi-analytische Strukturen auf Teilmengen $Z$ von $X$ an.
 
\Beg{Satz und Definition}\label{1.10}
   Sei $Z$ eine Teilmenge von $X$ und $\5A$ ein System lokal-analytischer Teilmengen von $X$ mit folgenden Eigenschaften:
   \BegEN
      \item\label{1.10.1}
         $Z = \bigcup_{A \in \5A} A$
      \item\label{1.10.2}
         $\forall A,B \in \5A$ ist $A \cap B$ sowohl $A$- als auch $B$-offen (dabei ist $A \cap B = \emptyset$ zugelassen)
   \EndEN
   Dann ist
   $$\5B := \{A' : A' \Text{ist $A$-offene Teilmenge eines} A \in \5A\}$$
   Basis einer quasi-analytischen Struktur $\5T$ auf $Z$. $\5A$ hei{\ss}t ein  
   \DEF{Erzeugendensystem}\index{Erzeugendensystem von $\5T$} von $\5T$ und wir bezeichnen
   mit $ {\rm Top}(\5A) := \5T$ die durch $\5A$ erzeugte Topologie.
\End

{\bf Beweis:} Man beachte: sind $A,B \in \5A$ und ist $C \subset A \cap B$, so gilt:
$$
   C \Text{ist $A$-offen} \ \iff \  C \Text{ist $A \cap B$-offen} \ \iff \  C \Text{ist $B$-offen}
$$
Daraus folgt, da{\ss} $A' \cap B' \in \5B$ f\"ur alle $A', B' \in \5B$. Also ist $\5B$ Basis einer Topologie $\5T$ auf $Z$. 
Offensichtlich ist $\5T$ eine quasi-analytische Struktur . \Qed

\Beg{Satz}\label{1.11}\parskip = 0pt
   Die Teilmenge $Z$ von $X$ besitze folgende Eigenschaft:
   
         Zu jedem $z \in Z$ gibt es eine offene $X$-Umgebung $U$ von $z$ derart, da{\ss} die $X$-Weg\-zu\-sam\-men\-hangs\-komponente
      $A$ von $Z \cap U$, die $z$ enth\"alt, eine analytische Teilmenge von $U$ ist.
   
   Dann erf\"ullt das System $\5A$ aller dieser lokal-analytischen Teilmengen $A$ von $X$ die Voraussetzungen 
   von \ref{1.10}  und  $\5T = {\rm Top}(\5A)$ ist eine Wege-vertr\"agliche quasi-analytische Struktur auf $Z$.
\End

{\bf Beweis:} Da{\ss} $\5A$ die Bedingung \ref{1.10}.\ref{1.10.1} erf\"ullt, folgt aus dem lokalen Wegzusammenhang lokal-analytischer Mengen. Wegen \ref{1.5}.\ref{1.5.3} ist $\5T$ Wege-vertr\"aglich.  \Qed

Wegen \ref{1.6} ist das System $\5A$ aus \ref{1.11} ein Erzeugendensystem f\"ur jede Wege-vertr\"agliche quasi-analytische Struktur auf $Z$.

Bezeichne $\5T$ die quasi-analytische Struktur auf $Y$ in \ref{1.4}.\ref{1.4.1}. Dann ist $\5A  = \{ \{z\} \times \C : z \in M\}$ ein
Erzeugendensystem von $\5T$.

Mit $\Theta$ bezeichnen wir die Garbe der holomorphen Vektorfelder auf $X$. F{\"u}r ein Vektorfeld $\theta$ auf $W \Opss X$
und $z \in W$ bezeichnen wir mit $\theta|_z$ den durch $\theta$ im Tangentialraum $\TT(X,z)$ definierten Tangentialvektor.

\Beg{Definition und Bemerkung}\label{1.12}
   Sei $(Z,\5T)$ eine Schichtung in $X$ und $z \in Z$. Dann ist f\"ur den komplexen Raum $(Z,\5T)$ der 
   \DEFI{Tangentialraum} 
   $\TT\big((Z,\5T),z\big) \subset \TT(X,z)$ in $z$ wohldefiniert.
\End

Ist in der Situation von \ref{1.12} $(A,U)$ ein $\5T$-Pl\"attchen mit $z \in A$, so ist $\TT\big((Z,\5T),z\big) = \TT(A,z)$ (zur Erinnerung an die Definition von $\TT(A,z)$ vgl.~den Beweis von \ref{1.14}).

\Beg{Definition}\label{1.13}
   Sei $(Z,\5T)$ eine Schichtung in $X$ und $W \Opss X$. Ein Vektorfeld $\theta \in \Theta(W)$ hei{\ss}t \DEF{tangentiell} 
   \index{tangentiell zu $(Z,\5T)$} oder
   \DEFI{parallel zu $(Z,\5T)$}, wenn gilt:
   $$\theta|_{z} \in \TT\big((Z,\5T),z\big) \quad \forall z \in W.$$
   Wir schreiben dann $\theta \| (Z,\5T)$.
\End

Wenn bei einer Betrachtung einer Schichtung $(Z,\5T)$ die Hervorhebung der quasi-analytischen Struktur nicht notwendig erscheint, lassen wir im Folgenden h\"aufig die Angabe von $\5T$ fort, sprechen also von der Schichtung $Z$ und notieren $\TT(Z,z)$,
$\theta \| Z$ und benutzen \"ahnliche Redeweisen bzw.~Notationen.

\Beg{Satz}\label{1.14}
   Sei $Z$ eine Schichtung in $X$, $W \Opss X$ und $\theta \in \Theta(W)$. Wenn es eine bez\"uglich der $Z$-Topologie dichte Teilmenge
   $M \subset Z \cap W$ gibt derart, da{\ss} $\theta|_{z} \in \TT(Z,z)$ f\"ur alle $ z \in M$, dann ist $\theta$ parallel zu $Z$.
\End

{\bf Beweis:} Es sei $z \in Z \cap W$ und $(A,U)$ ein Pl\"attchen von $Z$ mit $z \in A$. 
Wir d\"urfen annehmen, da{\ss} $U \subset W$ und da{\ss}  es Funktionen $f_{1},\ldots,f_{m} \in \5O(U)$ gibt, welche die Idealgarbe von $A$ \"uberall in  $U$ erzeugen. Dann ist der lineare Faserraum 
$$\TT(A) := \big\{(z,\eta) \in U \times \C^n : f_{\mu}(z)=0, d_{z}f_{\mu}(\eta |_{z})=0, \mu = 1,\ldots,m\big\}$$
der Faserraum der Tangentialr{\"a}ume an $A$. 
Nun ist die Behauptung offensichtlich. \Qed

\Beg{Bezeichnung}
   Es sei $Z$ eine Schichtung in $X$. 
   Dann notieren wir die $\5O_{X}$-Garbe der zu $Z$
   tangentiellen Vektorfelder mit $\Theta^Z$.
\End

Die Garbe $\Theta^Z$ ist involutiv:

\Beg{Satz}\label{1.16}
   Es sei $Z$ eine Schichtung in $X$, $W \Opss X$ und 
   $\theta,\eta \in \Theta^Z(W)$. Dann ist auch $\Clint \theta \eta \in \Theta^Z(W)$.
\End

{\bf Beweis:} Wir gehen von der Situation im Beweis von \ref{1.14} aus und sehen Vektorfelder als Derivationen an.
Mit $\5I$ notieren wir die Idealgarbe von $A$ in $U$. Dann gilt 
$\eta(f_{\mu}) \in \5I(U), \theta(f_{\mu}) \in \5I(U)$ f\"ur jedes $f_{\mu}$, also
$\Clint \theta \eta(f_{\mu}) = \theta\big(\eta(f_{\mu})\big) - \eta\big(\theta(f_{\mu})\big) \in \5I(U).$ \Qed

In Definition \ref{1.1} haben wir definiert, was wir unter einer quasi-analytischen Struktur auf einer Teilmenge $Z$ der Mannigfaltigkeit $X$ verstehen wollen. In der folgenden Definition verallgemeinern wir diese Definition in kanonischer Weise f{\"u}r Teilmengen  beliebiger komplexer R{\"a}ume:

\Beg{Definition}\label{1.17}
   Sei $Y$ ein komplexer Raum. Eine \DEFI{quasi-analytische Struktur} auf einer Teilmenge $Z$ von $Y$ ist eine
   Topologie $\5T$ auf $Z$, f\"ur die gilt:
   \\
   Zu jedem $z \in Z$ gibt es eine offene $\5T$-Umgebung $A$ von $z$ in $Z$ und eine offene $Y$-Umgebung 
   $U$ von $z$ in $Y$, derart da{\ss} $A$ eine analytische Teilmenge von $U$ ist und da{\ss} die $\5T$-Topologie 
   (vgl.~\ref{1.1})
   von $A$ mit der $Y$-Topologie von $A$ \"ubereinstimmt.
\End

Wir benutzen in der Situation von \ref{1.17} Bezeichnungen analog zu \ref{1.1}.

Im Weiteren benutzen wir die Begriffsbildung \ref{1.17} an einigen Stellen in der speziellen Situation, da{\ss} der komplexe Raum $Y$ eine Schichtung $(Z,\5T)$ von $X$ ist. 

\Beg{Definition}\label{1.18}
   Seien $(\Bul Z,\Bul{\5T})$ und $(Z,\5T)$ quasi-analytische Schichtungen von $X$. Dann hei{\ss}t 
   $(\Bul Z,\Bul{\5T})$ eine \DEFI{quasi-analytische Teilmenge} von  $(Z,\5T)$, wenn gilt:
   \BegEN
      \item $\Bul Z \subset Z$,
      \item $(\Bul Z,\Bul{\5T})$ ist eine quasi-analytische Schichtung von  $(Z,\5T)$.
   \EndEN
\End

\Beg{Satz}\label{1.19}
   Seien $(\Bul Z,\Bul{\5T})$ und $(Z,\5T)$ quasi-analytische Schichtungen von $X$. Dann gilt:
   $(\Bul Z,\Bul{\5T})$ ist genau dann eine quasi-analytische Teilmenge von  $(Z,\5T)$, 
   wenn es zu jedem $z \in \Bul Z$ ein $\Bul{\5T}$-Pl\"attchen $\Bul A$ und ein $\5T$-Pl\"attchen $A$ gibt mit
   $z \in \Bul A \subset A$.
\End

Der {\bf Beweis} ist klar. 

Ist in der Situation von \ref{1.18} $\Bul{Z'}$ eine Schicht von $\Bul Z$, so ist $\Bul{Z'}$ 
$\5T$-zusammenh\"angend, also in einer Schicht $Z'$ von $Z$ enthalten, genauer: quasi-analytische Teilmenge einer Schicht $Z'$ von $Z$.

\subsection{Beispiele}\label{Bsp}

Vorweg zwei Aussagen, welche Begr\"undungen bei der Diskussion einiger der nachfolgenden Beispiele liefern.

Die Standardstruktur einer lokal-analytischen Teilmenge ist Wege-vertr\"aglich, deshalb folgt mit \ref{1.6}:

\Beg{}\label{2.1}
   Sei $A$ eine lokal-analytische Teilmenge in $X$. Ist $\5T$ eine quasi-analytische Struktur auf $A$, 
   die von der Standardstruktur verschieden ist, so ist $\5T$ nicht Wege-vertr\"aglich.
\End

\Beg{Satz}\label{2.2}
   Die Teilmenge $Z$ von $X$ sei lokal $X$-wegzusammenh\"angend. Gibt es auf $Z$ eine Wege-vertr\"agliche 
   quasi-analytische Struktur $\5T$, so ist $Z$ eine lokal-analytische Teilmenge von $X$ und $\5T$ ist die Standardstruktur.
\End

{\bf Beweis:} Sei $z \in Z$ und $(A,U)$ ein $\5T$-Pl\"attchen zu $z$ gem\"a{\ss} \ref{1.5}.\ref{1.5.3}. Es existiert eine $X$-Umgebung
$W \subset U$ derart, da{\ss} $Z \cap W$ $X$-zusammenh\"angend ist. Es folgt: $A \cap W = Z \cap W$, $Z$ ist in $z$ analytisch.
\Qed

Ist die Teilmenge $Z$ in $X$ also lokal $X$-zusammenh\"angend und nicht lokal-analytisch, so kann sie keine Wege-vertr\"agliche quasi-analytische Struktur besitzen.

In den folgenden Beispielen gewinnen wir die quasi-analytischen Strukturen mit Hilfe von \ref{1.10}.

\BegRM{Beispiel}\label{2.3}
   Sei $X =Z =  \C$, $\5A := \big\{\{0\}, \C^*\big\}$. 
   Dann ist $\5T := {\rm Top}(\5A)$
   \index{Top}  eine von der Standardstruktur auf $\C$ verschiedene Struktur. 
   $\5T$ ist nicht Wege-vertr\"aglich.
\EndRM

\BegRM{}\label{2.4}
   Sei $X := \C^2$ und $$Z = \big\{(z,w) \in \C^2 : z \cdot w = 0\big\} = \big(\C \times \{0\}\big) \cup \big(\{0\} \times \C\big).$$
   Dann ist 
   $$\5T := {\rm Top}\big(\C \times \{0\}, \{0\} \times \C^*\big)$$
   eine von der Standardstruktur auf der analytischen Menge $Z$ verschiedene Struktur, insbesondere nicht Wege-vertr\"aglich. 
   $(Z,\5T)$ besteht aus zwei Schichten und ist eine Mannigfaltigkeit.
\EndRM

\BegRM{Beispiel}\label{2.5}
   Es sei $X := \C^2$ und
   $$Z_{0}:= \big(\C \times \{0\}\big) \cup \big(\{0\} \times \{w \in \C : |w-1|>1\}\big).$$
   $Z_{0}$ ist in $(0,0)$ nicht analytisch, aber lokal $X$-wegzusammenh\"angend. Deshalb kann $Z_{0}$ keine
   Wege-vertr\"agliche quasi-analytische Struktur tragen. Insbesondere ist die quasi-analytische Struktur
   $$\5T_{0} := {\rm Top} \big(\C \times \{0\}, \{0\} \times \{w \in \C : |w-1|>1\}\big)$$
   nicht Wege-vertr\"aglich. $(Z_{0},\5T_{0})$ besteht aus zwei Schichten und ist eine Mannigfaltigkeit.
\EndRM

Die Konstruktionsideen bei \ref{2.4} und \ref{2.5} kann man benutzen, um auch exotische zusammenh\"angende Schichtungen zu bekommen.

Im Folgenden bezeichnen wir $\C$ mit $\C_{z}$, wenn wir andeuten wollen, da{\ss} wir die Elemente von $\C$ mit $z$ bezeichnen; analog wird $\C_{z,w}$ definiert.

\BegRM{Beispiel}\label{2.6}
   Sei $X := \7P^2$. Dann ist
   $Z := \big\{ [z,w,t] \in \7P^2 : z \cdot w \cdot t = 0 \big\}$ eine komplexe Kurve im $\7P^2$ mit den irreduziblen Komponenten
   $$\7P^1_{1}:=\big\{[z,w,t] \in \7P^2 : z = 0 \big\},\quad
   \7P^1_{2}:=\big\{[z,w,t] \in \7P^2 : w = 0 \big\},\quad
   \7P^1_{3}:=\big\{[z,w,t] \in \7P^2 : t = 0 \big\}.$$
   $Z$ ist zusammenh\"angend, da
   $$[0,0,1] \in \7P^1_{1} \cap \7P^1_{2},\quad
     [0,1,0] \in \7P^1_{1} \cap \7P^1_{3},\quad  
     [1,0,0] \in \7P^1_{2} \cap \7P^1_{3}.$$

     Wir interpretieren wie \"ublich 
     $\C^2_{z,w}$ als
     Teilmenge von $\7P^2$:
     $$\C^2_{z,w}  = \big\{[z,w,t]  \in \7P^2 : t = 1\big\}.$$
     Damit ist $Z' := Z \cap \C^2_{z,w}$ die in Beispiel \ref{2.4} angegebene analytische Menge. 

      Sei $$W := \7P^2 \sm \big\{(z,w) \in \C^2_{z,w} : |z|^2 + |w|^2 \leq 1\big\}.$$
      Dann ist $$\5T := {\rm Top}\big(\C_{z} \times \{0\}, \{0\} \times \C^*_{w}, Z \cap W\big)$$
      eine von der Standardstruktur auf $Z$ verschiedene Struktur, insbesondere nicht Wege-vertr\"aglich.
      $(Z,\5T)$ besteht aus einer Schicht.
      
      Wir variieren nun die Menge $Z'$ wie in \ref{2.5}. Sei
      $$Z_{0}' := \big(\C_{z} \times \{0\}\big) \ \cup\  \big(\{0\} \times \{w \in \C_{w} : |w-1|>1\}\big).$$
      Wir ersetzen die Menge $Z'$ in $Z$ durch die Menge $Z_{0}'$ und erhalten so die Menge $Z_{0}$. Sie ist
      in $(0,0)$ nicht analytisch, aber lokal $X$-wegzusammenh\"angend. Deshalb kann $Z_{0}$ keine 
      Wege-vertr\"agliche quasi-analytische Struktur tragen. Insbesondere ist die quasi-analytische Struktur
      $$\5T_{0}:= {\rm Top}\big(\C_{z} \times \{0\}, \{0\} \times \{w \in \C_{w} : |w-1|>1\}, Z \cap W\big)$$
      auf $Z_{0}$ nicht Wege-vertr\"aglich. $(Z_{0},\5T_{0})$ besteht aus einer Schicht.
\EndRM

\BegRM{Beispiel (Kurve von Lissajous)}\label{2.8}
   Es sei $X := \C^2$ und 
   $$Z := \big\{(z,w) \in \C^2 : 4z^4 - 4z^2 + w^2 = 0\big\}$$
   Dann gilt:
   \BegEN
      \item Sing $Z = \{(0,0)\}$\label{L1}
      \item\label{L2} f\"ur $W := \{(z,w) \in \C^2 : |z| < 1\}$ ist $Z' := Z \cap W = Z_{+} \cup Z_{-}$, wobei
      $$Z_{\pm} := \{(z, \pm 2iz\sqrt{z^2-1}):z \in \C, |z|<1\};$$ $Z_{+}$ und $Z_{-}$ sind Untermannigfaltigkeiten
      von $W$, die sich in $(0,0)$ schneiden.
      \item\label{L3} f\"ur $\Mapping g \C {\C^2}$, $g(t) := (\cos t,\sin 2t)$ ist $g(\C) = Z$
      \item\label{L4} $Z \sm \{(0,0)\}$ ist zusammenh\"angend.
   \EndEN
   {\bf Beweis}:
   \\
   {\bf \ref{L1}}:
   Sei $f:=4z^4-4z^2+w^2 = 4z^2(z^2-1)+w^2$.    
    Es ist $df = (16z^3 - 8z)\,dz + 2w\,dw$, also
   \BegEA
    df|_{(z,w)} = 0 &\iff&8z(2z^2-1)=0 \Text{und} w=0\\
    &\iff& (z=0, w=0) \Text{oder} (z^2=1/2, w=0).
   \EndEA
   F\"ur $z^2 = 1/2, w=0$ ist $f(z,w) = -1 \not = 0.$ Also ist Sing $Z \subset \{(0,0)\}$.
   \\
   \ref{L1} und \ref{L2}: Sei $(z,w) \in W$. Dann gilt:
   \BegEA
      (z,w) \in Z &\iff& 4z^2(z^2-1) + w^2 = 0\\
      &\iff&(2z\sqrt{z^2-1}+iw)(2z\sqrt{z^2-1}-iw)=0\\
      &\iff&w =  2iz\sqrt{z^2-1} \Text{oder} w = - 2iz\sqrt{z^2-1}\\
      &\iff& (z,w) \in Z_{+} \cup Z_{-}.
   \EndEA
   \ref{L3}: Sei $t \in \C$ und $(z,w) = g(t)$. Dann gilt:
   \BegEA
      4z^2(z^2-1)+w^2 &=& 4(\cos^2t)(\cos^2t-1) + \sin^22t\\
      &=& -4\cos^2 t \cdot \sin^2 t + 4 \cos^2 t  \cdot \sin^2 t\\
      &=& 0
   \EndEA
   Folglich ist $g(\C) \subset Z$.
   
   Es sei $(z,w) \in Z$. Da $\Mapping \cos \C\C$ surjektiv ist (vgl.~\ref{7.7}),
   existiert ein $t \in \C$ mit $z = \cos t$. Dann folgt:
      $$w^2 = 4z^2(1-z^2) = 4\cos^2t\sin^2t = (\sin 2t)^2,$$
      also 
      $w = \sin 2t$ oder $w = - \sin 2t$, also $(z,w) = g(t)$ oder $(z,w) = g(-t)$. 
      
      \ref{L4}: Sei $\Z' := \{k \in \Z : k \Text{ist ungerade$\}$}$. Dann gilt:
      $g(t) = (0,0) \iff t=k\pi/2$, $k \in \Z'$. Also ist $Z \sm \{(0,0)\} = g\big(\C \sm (\Z'\cdot \pi/2)\big).$
\EndRM

\BegRM{Fortsetzung von \ref{2.8}}\label{2.9}
   $Z$ ist eine analytische Menge mit einer Singularit\"at in $(0,0)$. Dann ist
   $$\5T := {\rm Top}\big(Z \sm \{(0,0)\}, Z_{+}, Z_{-} \sm \{(0,0)\}\big)$$
   eine von der Standardstruktur verschiedene quasi-analytische Struktur auf $Z$,
   insbesondere nicht Wege-vertr\"aglich. $(Z,\5T)$ ist eine Schicht und zwar eine Riemannsche Fl\"ache.
   
   Wir variieren nun die Menge $Z$ analog wie in \ref{2.5}. Sei
   $$Z'_{0} := Z_{+} \cup \Big(Z_{-} \cap \big\{(z,w) \in \C^2 : |z-1/3| > 1/3\big\}\Big).$$
   Wir ersetzen die Menge $Z'$ in $Z$ durch die Menge $Z_{0}'$ und erhalten so die Menge $Z_{0}$. Sie ist in 
   $(0,0)$ nicht analytisch, aber lokal $X$-wegzusammenh\"angend, 
   kann also keine Wege-vertr\"agliche quasi-analytische Struktur tragen.
   Insbesondere ist die quasi-analytische Struktur
   $$\5T_{0} := {\rm Top}\big(Z \cap \{(z,w) \in \C^2 : |z-1/3| > 1/3\}, Z_{+}, Z_{-} \cap 
   \{(z,w) \in \C^2 : |z-1/3| > 1/3\}\big)$$
   auf $Z_{0}$ nicht Wege-vertr\"aglich. $(Z_{0},\5T_{0})$ ist eine Schicht und zwar eine Riemannsche Fl\"ache.
\EndRM

Im folgenden Beispiel benutzen wir \ref{1.2}.\ref{1.2.2} zur Konstruktion einer quasi-analytischen Struktur.

\BegRM{Beispiel}\label{2.10}
     Wir betrachten den eindimensionalen Torus $T := \C/\Z^2$ sowie den zweidimensionalen Torus
   $X := \C^2/\Z^4 = T \times T$. Die nat\"urliche Projektion $\Mapping p {\C^2}X$ ist lokal biholomorph. Sei $a \in \R$ irrational
   und $\Mapping \phi \C X$, $\phi(z) := p(z,az)$. 
   Offensichtlich ist $\phi$ eine holomorphe Immersion. $\phi$ ist auch injektiv:
   \BegEA
      \phi(z) = \phi (w) & \iff & p(z,az) = p(w,aw)\\
      &\iff& z-w =: k \in \Z^2,  ak \in \Z^2\\
      &\iff& k=0 \\
      &\iff& z=w.
   \EndEA

   Sei $\5T$ die durch $\phi$ gem\"a{\ss} \ref{1.2}.\ref{1.2.2} definierte quasi-analytische Struktur auf $Z = \phi(\C)$.
   Dann ist $(Z,\5T)$ eine Schicht.

      $Z$ ist dicht in $X$ (bei der folgenden Argumentation gehen wir reell vor und notieren f\"ur komplexe Zahlen 
      $z \in \C$ die reellen Komponenten mit $z_{1},z_{2}$):
      Sei $(u,v) \in \C^2$, $\eps > 0$. Wegen \ref{7.5} gibt es ein $z = (z_{1},z_{2}) \in \C$ und Zahlen $k_{1}, k_{2}, 
      l_{1}, l_{2} \in \Z$ mit:
      $$z_{1} = u_{1} + k_{1}, \quad |az_{1} - (v_{1} + l_{1})| < \eps;$$
      $$z_{2} = u_{2} + k_{2}, \quad |az_{2} - (v_{2} + l_{2})| < \eps.$$
      Man schlie{\ss}t wie in  \ref{7.6}: $Z$ ist dicht in $X$.
      
      Die quasi-analytische Struktur  ist sehr ordentlich, insbesondere  Wege-vertr\"aglich:
      
      Sei $(\eta,\xi) \in \TT \times \TT$ und $\eps \in \R$ mit $0 < \eps < 1/2$. 
      Wir betrachten einen Punkt $(u,v) \in \C^2$
      mit $p(u,v)= (\eta,\xi)$. Sei (wieder in reeller Notation)
      $$Q_{\eps} := \Opint {u_{1}-\eps}{u_{1}+\eps} \times \Opint{u_{2}-\eps}{u_{2}+\eps}
      \times  \Opint {v_{1}-\eps}{v_{1}+\eps} \times \Opint{v_{2}-\eps}{v_{2}+\eps} \Opss \C^2$$
      und $W_{\eps} := p(Q_{\eps}) \Opss X$. Die Abbildung $\Mapping p {Q_{\eps}}{W_{\eps}}$ ist biholomorph und
      $$p^{-1}(W_{\eps}) = \bigcup_{k,l \in Z^2} \big(Q_{\eps}+ (k,l)\big)$$
      (dabei sind die Mengen in der Vereinigung paarweise disjunkt).
      
      Sei $G := \{(z,az) : z \in \C\}$. F\"ur alle $k,l \in \Z^2$ ist $G_{\eps,k,l}:= G \cap \big(Q_{\eps} + (k,l)\big)$
      eine $\big(Q_{\eps} + (k,l)\big)$-abgeschlossene und konvexe Teilmenge von $Q_{\eps} + (k,l)$. Sei
      $$A_{\eps,k,l} := p\big(Q_{\eps} +(k,l)\big).$$
      Weil $\phi$ injektiv ist, sind die $A_{\eps,k,l}$ paarweise disjunkt. Sie bilden die $\5T$-Komponenten von 
      $Z \cap W_{\eps}$. Jedes $A_{\eps,k,l}$ ist $W_{\eps}$-abgeschlossen.
\EndRM

\subsection{Analytischer Inhalt}\label{AI}

In diesem Abschnitt wollen wir f\"ur eine Teilmenge von $X$ die N\"ahe zur Analytizit\"at beschreiben und f\"ur geeignete Teilmengen eine nat\"urliche quasi-analytische Struktur einf\"uhren.

Im Weiteren sei $Z$ eine Teilmenge von $X$.

\Beg{Definition}
   Die Menge $\5I(Z)$ aller lokal-analytischen Teilmengen $A$ von $X$ mit $A \subset Z$ nennen wir den
   \DEF{analytischen Inhalt}\index{analytischer Inhalt}\index{Inhalt!analytischer} von $Z$. Ferner nennen wir
   $$\dim Z := \max\{\dim A : A \in \5I(Z)\}$$
   die \DEFI{analytische Dimension}\index{Dimension!analytische} von $Z$.
\End

\Beg{}
   Sei $A$ eine lokal-analytische Teilmenge von $X$. Dann stimmen f\"ur $A$ die analytische Dimension
   und die \"ubliche (komplexe) Dimension \"uberein.
\End

\Beg{}\label{3.3}
   Sei $\5T$ eine quasi-analytische Struktur auf $Z$. Dann gilt:
   \BegEN
      \item Ist $(A,U)$ ein $\5T$-Pl\"attchen, so ist $A \in \5I(Z)$.
      \item $\dim(Z,\5T) \leq \dim Z$.\label{3.3.2}
   \EndEN
\End

\Beg{Satz}\label{3.4}
   Sei $\5T$ eine abz\"ahlbare quasi-analytische Struktur auf $Z$. Dann ist $\dim(Z,\5T) = \dim Z$.
\End

{\bf Beweis:} Wir f\"uhren die Annahme $p := \dim(Z,\5T) < q := \dim Z$ zum Widerspruch.

Es existiert eine zusammenh\"angende $q$-dimensionle Untermannigfaltigkeit $B$ von $X$ mit $B \subset Z$. Weil $(Z,\5T)$ ein Lindel\"of-Raum ist, gibt es eine abz\"ahlbares System $\5A$ von $\5T$-Pl\"attchen $(A,U)$ mit 
$Z \subset \bigcup_{A \in \5A} A$. F\"ur jedes $A \in \5A$ ist $B \cap A$ eine niederdimensionale lokal-analytische Teilmenge von $B$, also bzgl.~$B$ eine Lebesguesche Nullmenge. Weil $\5A$ abz\"ahlbar ist, mu{\ss} $B$ bzgl.~$B$ eine Lebesquesche Nullmenge sein. Widerspruch! \Qed

Da{\ss} die Forderung der Abz\"ahlbarkeit in \ref{3.4} wesentlich ist, zeigen die Beispiele \ref{1.4}. 
Genauer gilt: es bezeichne $\5T$ jeweils die dort angegebene quasi-analytische Struktur. Dann gilt:

\Beg{}
   Ist $M = \C$ in \ref{1.4}.\ref{1.4.1}, so ist $\dim(\C^2,\5T) = 1$ und $\dim\C^2=2$.
   --- Ist 
   $Y = X$ in \ref{1.4}.\ref{1.4.2}, so ist $\dim(X,\5T) = 0$ und $\dim X = n$.
\End

F\"ur die Extremf\"alle $\dim Z = n$ bzw.~$\dim Z = 0$ der analytischen Dimension gilt:

\Beg{}\label{3.6}\
   \BegEN
      \item $\dim Z = n \iff Z$ besitzt $X$-innere Punkte.
      \item $\dim Z = 0 \iff $ die diskrete Topologie ist die einzige quasi-analytische Struktur auf $Z$.
   \EndEN
\End

\BegRM{Beispiel}\label{3.7}
   Sei $X:=\C^2$ und $Z:=\{(u,v) : u,v \in \R\}$. Weil die reell 2-dimensionale reell-analytische Fl\"ache $Z$
   keine lokal-analytische Menge der (komplexen) Dimension 1 enth\"alt, gilt $\dim Z = 0$.
\EndRM

\Beg{}\label{3.8}
   Sei $\5T$ eine quasi-analytische Struktur auf $Z$. 
   \BegEN
      \item\label{3.8.1} Ist $\dim(Z,\5T) = n$, so enth\"alt $Z$ $X$-innere Punkte.
      \item\label{3.8.2} Ist $\5T$ abz\"ahlbar und enth\"alt $Z$ $X$-innere Punkte, so ist $\dim(Z,\5T) = n$.
   \EndEN
\End

{\bf Beweis:} Ad \ref{3.8.1}: Sei $(A,U)$ ein $\5T$-Pl\"attchen von $Z$ der Dimension $n$. Dann besitzt $A$ $X$-innere Punkte.
Ad \ref{3.8.2}: Wegen \ref{3.4} ist $\dim(Z,\5T) =  \dim Z$. \Qed

Aus \ref{3.4} folgt ferner:

\Beg{}\label{3.9}
   Seien $\5T, \Bul{\5T}$ abz\"ahlbare quasi-analytische Strukturen auf $Z$. 
   Dann ist $\dim(Z,\5T) = \dim(Z,\Bul{\5T})$.
\End

Wir nehmen jetzt an, da{\ss} die Situation von \ref{3.9} vorliegt. Dann sei $\tilde Z$ die Menge aller $z \in Z$ mit folgender Eigenschaft:

   \quad {\sl Es existiert eine offene $\5T$-Umgebung $V$ von $z$, die auch $\Bul{\5T}$-offen ist, und f{\"u}r die 
   $\5T|_{V}=\Bul{\5T}|_{V}$ gilt.}

Diese Forderung ist \"aquivalent zur folgenden Forderung:

\quad {\sl Es existiert ein $\5T$-Pl\"attchen $(A,U)$ mit $z \in A$, welches gleichzeitig ein $\Bul{\5T}$-Pl\"attchen ist.}

\Beg{}\label{3.10}
   Die Menge $\tilde Z$ ist $\5T$-offen und $\Bul{\5T}$-offen sowie $X$-dicht in $Z$.
\End

{\bf Beweis:}
Die Offenheitsaussage ist klar. Zum Nachweis der Dichtheit sei $z \in Z$, $W$ eine beliebige $X$-offene Umgebung von $z$ und $q := \dim(Z \cap W)$. Wir betrachten ein $\5T$-Pl\"attchen $(A,U)$ mit $U \subset W$. Es sei $B$ eine Zusammenhangskomponente von $A \sm \Sing A$ der Dimension $q$. Mit dem Argument aus dem Beweis von \ref{3.4} gibt es ein $\Bul{\5T}$-Pl\"attchen $(\Bul A, \Bul U)$ und ein $z' \in B \cap \Bul A$ mit 
$B_{z'} \subset (\Bul A)_{z'}$. Dann mu{\ss} $\dim{ \Bul A}_{z'} = q$ und $B_{z'}$ eine irreduzible Komponente von 
$(\Bul A)_{z'}$ sein. Indem wir gegebenenfalls $z'$ gegen einen Nachbarpunkt aus $B$ austauschen, d\"urfen wir 
$B_{z'} = (\Bul A)_{z'}$ und dann sofort $\Bul A$ als $B$-offene Teilmenge von $B$ annehmen. 
Also ist $z' \in \tilde Z$.\Qed

Eine einfache Modifikation des Beweises von \ref{3.10} liefert

\Beg{Satz}\label{3.11}
   Sei $p = \dim Z$ und $Z_{p}$ die Vereinigung der irreduziblen Komponenten von $(Z,\5T)$ der Dimension $p$.
   Dann ist $\tilde Z \cap Z_{p}$ $\5T$-dicht in $Z_{p}$.
\End

Wie unsere Beispiele \ref{2.3}, \ref{2.4}, \ref{2.6} und \ref{2.9} zeigen, ist unser Ergebnis \ref{3.10} bzw.~\ref{3.11}
die bestm\"ogliche Eindeutigkeitsaussage f\"ur abz\"ahlbare quasi-analytische Strukturen.

Im Folgenden wollen wir mit Hilfe von $\5I(Z)$ und Satz \ref{1.10} f\"ur geeignete Mengen $Z$ eine nat\"urliche quasi-analytische Struktur auf $Z$ konstruieren. 
Um \ref{1.10}.\ref{1.10.2} zu erf\"ullen, m\"ussen wir eine Maximalit\"atsforderung stellen. Dazu m\"ussen wir gewisse Teilmengen von $\5I(Z)$ betrachten. Zun\"achst einige Bezeichnungen:

\Beg{Bezeichnung}
   Sei $Y$ ein komplexer Raum. Dann bezeichne
   $$\dimMin Y := \min\{\dim_{y}Y : y \in Y\}$$
   die \DEFI{Minimaldimension} von $Y$.
\End

\Beg{Bezeichnung}
   F\"ur $p \in \N_{0}$ sei
   $$\5I_{p}(Z) := \{ A \in \5I(Z) : \dimMin A \geq p\}.$$
\End

Nat\"urlich ist $\5I_{p}(Z) = \emptyset$ f\"ur $p > \dim Z$.

Nun zur zentralen Definition (dabei bezeichnen wir, wie bisher schon, mit $A_z$ den von $A$ in $z$ definierten Mengenkeim):

\Beg{Definition}\label{3.14}
   F\"ur $p \in \N_{0}$ sei 
   $$\5A_{p} =\5A_{p}(Z)\  :=\  \big\{ A \in \5I_{p}(Z) : \Text{f\"ur alle $B \in \5I_p(Z)$ gilt:} [B_{z} \subset A_{z} \ \forall
   z \in B \cap A]\big\}$$
   Ist $A \in \5A_{p}$ zusammenh\"angend (was man durch geeignete Verkleinerung stets erreichen kann), so hei{\ss}t $A$ eine \DEF{$p$-dicke $Z$-Scheibe}\index{p-dicke Scheibe@$p$-dicke Scheibe}\index{dick}\index{Scheibe!dicke}
\End

\BegRM{Beispiel}\label{3.15}
   Sei $Z$ eine lokal-analytische Teilmenge von $X$ und $p := \dimMin Z$. Dann gilt:
   $$\5A_{0}(Z) = \5A_{1}(Z) = \ldots = \5A_{p}(Z) = \{A : A \Text{ist $Z$-offene Teilmenge von} Z\}$$
\EndRM

Ist in \ref{3.14} auch $B \in \5A_{p}$, so gilt f\"ur alle $z \in A \cap B$, dass $B_{z} = A_{z}$. Daraus folgt:

\Beg{Definition und Satz}\label{3.16}
   $Z$ hei{\ss}t eine \DEF{$p$-dicke schwach-analytische Teilmenge}
   \index{p-dicke schwach analytische Teilmenge@$p$-dicke schwach-analytische Teilmenge}
    von $X$, wenn $Z = \bigcup_{A \in \5A_{p}} A$.
---
In diesem Fall erf\"ullt $\5A_{p}$ die Voraussetzungen von \ref{1.10}, definiert also eine quasi-analytische Struktur 
$\5T_{p}$ auf $Z$, die 
\DEF{$p$-Standardstruktur}\index{p-Standardstruktur@$p$-Standardstruktur}\index{Standardstruktur}.
\End

\Beg{}
  Sei $Z$ eine lokal-analytische Teilmenge von  $X$ und $p := \dimMin Z$. F\"ur jedes $q \leq p$ ist dann $Z$ eine 
  $q$-dicke schwach-analytische Teilmenge von $X$ und $\5T_{q}$ ist die \"ubliche Standardtopologie.
\End

Zum {\bf Beweis} siehe \ref{3.15}. 

\Beg{Satz}\label{3.18}
   Sei $\5T$ eine Wege-vertr\"agliche quasi-analytische Struktur auf $Z$. Dann ist $Z$ eine $0$-dicke schwach-analytische Teilmenge von $X$ und $\5T = \5T_{0}$.
\End

{\bf Beweis:} Bezeichne $\5A$ das in \ref{1.11} genannte Erzeugendensystem von $\5T$. Weil lokal-analytische Mengen lokal $X$-wegzusammenh\"angend sind, folgert man sofort dass $\5A \subset \5A_{0}(Z)$. \Qed

Wie Beispiel \ref{3.15} zeigt, ist die Zahl $p$ bei einer $p$-dicken schwach-analytischen Teilmenge $Z$ von $X$ nicht eindeutig bestimmt. Allerdings gilt:

\Beg{Satz}\label{3.19}
   Sei $Z$ eine zugleich $p$-dicke und $q$-dicke schwach-analytische Teilmenge von $X$. 
   Dann ist $\5A_{p}=\5A_{q}$ und folglich $\5T_{p}=\5T_{q}$.
\End

Zum Beweis von \ref{3.19} zeigen wir zun\"achst:

\Beg{}\label{3.20}
   Sei $Z$ eine $p$-dicke schwach-analytische Teilmenge von $X$ und $q \leq p$. 
   Dann ist $\5A_{q} \subset \5A_{p}.$
\End

{\bf Beweis}: Sei $B \in \calA_{q}$. Wir m\"ussen zeigen: $\dimMin B \geq p$. Wir nehmen an: $\dim B < p$. Dann gibt es einen regul\"aren Punkt $x$ von $B$ mit $\dim_{x}B < p$. Sei $A \in \calA_{p}$ mit $x \in A$. Dann ist 
$A_{x}\subset B_{x}$, insbesondere $\dim_xB \geq p$. Widerspruch!

{\bf Beweis} von \ref{3.19}: Wir d\"urfen annehmen dass $q \leq p$. Wegen \ref{3.20} ist $\calA_{q} \subset \calA_{p}$. Sei 
$A \in \calA_{p}$ und $x \in A$. Es existiert ein $B \in \calA_{q}$ mit $x \in B$. Weil auch $B \in \calA_{p}$ gilt, 
ist $B_{z} = A_{z}$, und es existiert eine offene $X$-Umgebung $U$ von $z$ mit $A \cap U = B \cap U$. Aus dieser \"Uberlegung folgert man da{\ss} f\"ur alle $C \in \5I_{q}$ gelten mu{\ss}: $C_{z} \subset A_{z}\ \forall z \in C \cap A$.
Also ist $A \in \5A_{q}$. \Qed

\Beg{Definition}
   $Z$ hei{\ss}t eine \DEFI{schwach-analytische Teilmenge} von $X$, wenn es ein $p \in \N_{0}$ gibt derart,
   da{\ss} $Z$ eine $p$-dicke schwach-analytische Teilmenge von $X$ ist. Die $p$-Standardstruktur von $Z$,
   die im Sinne von \ref{3.19} nicht von $p$ abh\"angt, hei{\ss}t kurz \DEFI{Standardstruktur}. 
   Das kleinste zul\"assige $p$ nennen wir die \DEFI{Feinheit}  von $Z$ und schreiben daf\"ur
   ${\rm fein}(Z)$\index{fein(Z)@fein$(Z)$}.
\End

Je kleiner fein$(Z)$ ist, desto sch\"arfer ist  in \ref{3.14} die Bedingung $B_{z} \subset A_{z} \ \forall z \in B \cap A$.

\Beg{}\label{3.22}
   Besitzt $Z$ eine Wege-vertr\"agliche quasi-analytische Struktur, so ist $Z$ schwach-analytisch von der Feinheit 0.
\End

Zum {\bf Beweis} vgl.~\ref{3.18}.

Insbesondere ist die Menge $Z$ aus \ref{2.10} schwach-analytisch von der Feinheit 0.

Man sieht sehr leicht:

\Beg{}\label{3.23}
   Die Mengen $Z_{0}$ in \ref{2.5}, \ref{2.6} und \ref{2.9} sind schwach-analytisch von der Feinheit 0. Die dort
   angegebenen quasi-analytischen Strukturen sind die Standardstrukturen.
\End

\BegRM{Beispiel}\label{3.24}
   Sei $X := \C^3$ und 
   $Z:=\big(\C^2_{z,w} \times \{0\}\big) \ \cup \ \big(\{0\} \times \{(w,t) \in \C^2_{t,w}:|t|>|w|\}\big).$ 
   
   Sei $\5A := \big\{\C^2_{z,w} \times \{0\}\big\},\{0\} \times \{(w,t) \in \C^2_{w,t} : |t| > |w|\}\big\}$.
   Weil $Z$ die komplexe Gerade $\{(0,0)\} \times \C_{t}$ enth\"alt, gilt: 
   $(0,0,0) \not\in \bigcup_{A \in \5A_{0}}A =\bigcup_{A \in \5A_{1}}A$. Es gilt aber $\5A \subset \5A_{2}$. Also ist
   $Z$ eine schwach-analytische Teilmenge von $X$ von der Feinheit 2.
\EndRM

\BegRM{Beispiel}
   Sei $X := \C^2$ und $Z := \bigcup_{\alpha \in \Z} G_{\alpha}$, wobei 
   $G_{\alpha} := \{(z,w) \in \C^2 : \alpha z - w = 0\}$ ist. Weil durch den Punkt $(0,0)$ unendlich viele Geraden
   verlaufen, ist $(0,0) \not\in \bigcup_{A \in \5A_{p}}A \quad \forall p \in \N_{0}$. Allerdings definiert etwa
   $$\5A := \{G_{0}\} \cup \Big\{G_{\alpha} \sm \{(0,0)\} : 0 \not= \alpha \in \Z\Big\}$$ gem\"a{\ss} \ref{1.10}
   eine quasi-analytische Struktur auf $Z$.
\EndRM

\Beg{Satz}\label{3.26}
   Sei $Z$ eine schwach-analytische Teilmenge von $X$ und bezeichne $\5T$ die Standardstruktur. Dann ist
   $\dim(Z,\5T) = \dim Z$.
\End

{\bf Beweis:} Sei $q:= \dim(Z,\5T)$, $r := {\rm fein}(Z)$ und $p := \dim Z$. Wegen \ref{3.3}\ref{3.3.2}
 ist
$q \leq p$. Wir nehmen an, da{\ss} $q < p$. Dann sei $B \in \5I_{p}(Z)$ eine zusammenh\"angende Untermannigfaltigkeit von $X$ der Dimension $p$. Sei $z \in B$. Es existiert ein $A \in \5A_{r}(Z)$ mit $z \in A$. Nun gilt: $r \leq q < p$. Also ist $B \in \5I_{r}(Z)$, und es folgt: $B_{z} \subset A_{z}$, $\dim_{z}A \geq p$. Widerspruch!
\Qed

\subsection{Quasi-analytische Zerlegungen}\label{QaZ}

\Beg{Definition}
   Eine quasi-analytische Schichtung $(X,\5T)$ auf $Z = X$ nennen wir auch eine 
   \DEFI{quasi-analytische Zerlegung}\index{Zerlegung!quasi-analytische} von $X$. 
   Die zugeh\"origen Schichten, also die
   $\5T$-Zusammenhangskom\-po\-nen\-ten von $X$ mit der durch $\5T$ induzierten Topologie, 
   hei{\ss}en \DEF{Bl\"atter}\index{Blatt} der Zerlegung. Die Menge
   $\5D$ der Bl\"atter hei{\ss}t der \DEFI{Bl\"atterraum} der Zerlegung.
   Mit $\5D(x)$ bezeichnen wir dasjenige Blatt von $\5D$, welches $x$ enth{\"a}lt.
\End

\Beg{}\label{4.2}\
   \BegEN
      \item
         Auf $X$ liege eine quasi-analytische Zerlegung vor. Dann ist $X$ die disjunkte Vereinigung aller Bl\"atter
         dieser Zerlegung.
      \item
         Sei umgekehrt $X$ disjunkte Vereinigung von Schichten. Dann bildet die Menge $\5D$ aller dieser Schichten
         den Bl\"atterraum einer quasi-analytischen Zerlegung von $X$.
   \EndEN
\End

Der Beweis ist klar.

Wegen \ref{4.2} notieren wir im Folgenden eine quasi-analytische Zerlegung durch Angabe des
Bl\"atterraumes $\5D$ und bezeichnen die zugeh{\"o}rige Schichtung mit $(X,\5D)$.

Die Dimension des komplexen Raumes $(X,\5D)$ bezeichnen wir mit $\dim \5D$. Entsprechend ist $\dimMin \5D$
zu verstehen. Die Garbe der zu $(X,\5D)$ tangentiellen Vektorfelder bezeichnen wir mit 
$\Theta^{\5D}$\index{ThetaD@$\Theta^{\5D}$}

Es sei $U \Opss X$ zusammenh\"angend. Dann induziert $\5D$ eine quasi-analytische Struktur auf $U$. Wir notieren sie mit $\5D|_{U}$ und nennen sie die \DEFI{Beschr\"ankung} von $\5D$ auf $U$. 
Die Bl\"atter von $\5D|_{U}$ sind die $\5D$-Zusammenhangskomponenten der Mengen $A \cap U$, $A \in \5D$.

\Beg{Definition}\label{4.3}
   Eine quasi-analytische Zerlegung $\5D$ auf $X$ hei{\ss}e
   \BegIT
      \item
         \DEF{reindimensional}\index{Zerlegung!reindimensionale}, wenn $\dim \5D = \dimMin \5D$ ist, d.h.~wenn
         alle $A \in \5D$ reindimensionale komplexe R\"aume der gleichen Dimension sind
      \item
         \DEF{glatt}\index{Zerlegung!glatte}, wenn alle $A \in \5D$ Mannigfaltigkeiten sind,
      \item
         \DEF{(lokal)-analytisch}\index{Zerlegung!(lokal)-analytische}, wenn alle $A \in \5D$ (lokal)-analytische
         Teilmengen von $X$ mit der Standardstruktur sind
      \item
         \DEF{schwach-analytisch}\index{Zerlegung!schwach-analytische}, wenn es ein $p \in \N$ gibt, so da{\ss}
         jedes $A \in \5D$ eine $p$-dicke schwach-analytische Teilmenge von $X$ ist - das kleinste derartige $p$
         hei{\ss}t die \DEF{Feinheit}\index{Feinheit!einer Zerlegung} von $\5D$ und wird mit 
         {\rm fein}$(\5D)$\index{fein(D)@fein$(\5D)$} notiert
      \item\DEF{Wege-vertr\"aglich}\index{Zerlegung!Wege-vertr\"agliche} 
      bzw.~\DEF{(sehr) ordentlich}\index{Zerlegung!(sehr) ordentliche}, wenn alle $A \in \5D$ die entsprechende
      Eigenschaft haben.
   \EndIT
\End

Zur Erinnerung:
$$\Text{ordentlich} \ArrowComment{\IfThen}{\ref{1.9}}\Text{Wege-vertr\"aglich}
\ArrowComment{\IfThen}{\ref{3.22}}\Text{schwach-analytisch von der Feinheit 0}$$

Sei $(Z,\5T)$ eine Schichtung von $X$. Dann gilt (vgl.~Beweis von \ref{1.3}.\ref{1.3.2}): $Z$ ist genau dann eine analytische Teilmenge von $X$ und $\5T$ die Standardstruktur, wenn die Inklusion $\Mapping \iota Z X$ bez\"uglich $\5T$ eigentlich ist.

Deshalb sagen wir auch in Fortsetzung von \ref{4.3}:

\Beg{Definition} Eine quasi-analytische Zerlegung $\5D$ hei{\ss}t
\BegIT
   \item
      \DEF{eigentlich}\index{Zerlegung!eigentliche}, wenn $\5D$ analytisch ist;
   \item
      \DEF{lokal eigentlich}\index{Zerlegung!lokal eigentliche} wenn es zu jedem $x \in X$ eine offene
      zusammenh\"angende $X$-Umgebung $U$ von $x$ gibt, derart da{\ss} $\5D|_{U}$ eigentlich (also analytisch) ist.
\EndIT
\End

Ist $M = \C$ in \ref{1.4}.\ref{1.4.1}, so erhalben wir ein sehr einfaches Beispiel einer quasi-analytischen Zerlegung von $\C^2$. Sie ist reindimensional, glatt, analytisch und sehr ordentlich.
--- Wir verallgemeinern das Beispiel:

\Beg{}\label{4.5}
   Es seien $D_{1} \Opss \C^p$ und $D_{2} \Opss \C^q$ zwei Gebiete (d.h.~offen und zusammenh\"angend), 
   ferner sei
   $D := D_{1} \times D_{2}$. Dann ist $\5D := \{D_{1}\times \{w\} : w \in D_{2}\}$ eine quasi-analytische Zerlegung 
   von $D$. Sie ist rein $p$-dimensional, glatt, analytisch und sehr ordentlich.
\End

\Beg{}\label{4.6}
   Sei $\Mapping f X W$ eine holomorphe Abbildung in den komplexen Raum $W$. Dann sind alle Niveaumengen
   von $f$ (das sind die Zusammenhangskomponenten [bez{\"u}glich der $X$-Topologie] der Fasern $f^{-1}(w)$, $w \in W$)
   analytische Teilmengen von $X$ und wegen \ref{1.3} auf kanonische Weise Schichten.
      In diesem Sinne  ist
   $$\5D_f:=\{A : A \Text{\sl ist eine Niveaumenge von} f\}$$ 
   eine  quasi-analytische Zerlegung von $X$. Sie ist analytisch und sehr    ordentlich.
\End

Ist die Abbildung in \ref{4.6} offen, so d\"urfen wir $f$ als surjektiv voraussetzen. Da ein komplexer Raum $Y$ genau dann irreduzibel ist, wenn $Y \sm \Sing Y$ zusammenh\"angend ist (vgl.~\cite{KK}), folgert man, da{\ss} $W$ irreduzibel und lokal-irreduzibel ist. Indem wir $W$ gegebenenfalls durch seine Normalisierung ersetzen, d\"urfen wir $W$ sogar als normal voraussetzen. Es folgt, da{\ss} $\5D$ rein $p$-dimensional ist, wobei $p := n - \dim W$ 
(vgl.~\cite[49.16]{KK}).

Sei $\5F$ eine $p$-dimensionale regul\"are holomorphe Bl\"atterung auf $X$.\footnote{Hier und im Folgenden verweisen wir f\"ur Grundbegriffe der Theorie der holomorphen Bl\"atterungen auf \cite{Reiffen1}}
$\5F$ ist definiert durch eine regul\"are involutive analytische Untergarbe $\Theta_{\5F}$ von $\Theta$ vom Rang
$p$ bzw.~eine regul\"are involutive analytische Untergarbe $\Omega_{\5F}$ von $\Omega$ vom Rang $q:=n-p$ (dabei sei $\Omega$ die Garbe der holomorphen Pfaffschen Formen auf $X$). Die Garben $\Theta_{\5F}$ und 
$\Omega_{\5F}$ sind dual zueinander. $\5F$ kann lokal definiert werden durch gewisse holomorphe Submersionen, d.h.~zu jedem $x \in X$ existiert eine offene zusammenh\"angende Umgebung $U$ von $x$ und eine holomorphe Submersion $\Mapping f U W$ auf eine komplexe Mannigfaltigkeit mit $\Omega_{\5F}|_{U} = f^*(\Omega_{W})$.
Die durch die lokalen holomorphen Submersionen $\Mapping f U W$ gem\"a{\ss} \ref{4.6} definierten quasi-analytischen Zerlegungen auf den Gebieten $U$ in $X$ k\"onnen zu einer eindeutig definierten quasi-analytischen Zerlegung $\5D_{\5F}$\index{DF@$\5D_{\5F}$} auf $X$ verklebt werden. Lokal hat $\5D_{\5F}$ die Form von Beispiel \ref{4.5}.

\Beg{}\label{4.7}
   Sei $\5F$ eine $p$-dimensionale regul\"are holomorphe Bl\"atterung auf $X$. Dann ist $\5D_{\5F}$ glatt, 
   rein $p$-dimensional und sehr ordentlich.
\End

K.~Spallek hat gezeigt (vgl.~\cite{Spallek}):

\Beg{Satz}\label{4.8}
   Sei $\Theta'$ eine involutive koh\"arente Untergarbe von $\Theta$. Dann gibt es genau eine glatte 
   quasi-analytische Zerlegung $\5D'$ von $X$ mit der Eigenschaft:
   $$\Theta'|_{x} = \4T\big(x,\5D'(x)\big) \quad \forall x \in X.$$
\End
Wir nennen diese eindeutig bestimmte Zerlegung $\5D'$ die \DEFI{Spallek-Zerlegung} zu $\Theta'$.

\Beg{Satz}\label{4.9}
   In der Situation von  \ref{4.8} ist die Spallek-Zerlegung $\5D'$ von $\Theta'$ sehr ordentlich, insbesondere Wege-vertr{\"a}glich.
\End

{\bf Beweis:}
Es sei $A \in \5D'$, $\dim A = p$ und $x_{0} \in A$. Im Fall von $p = \rank \Theta'$ liegt in einer offenen $X$-Umgebung von $x_{0}$ die Situation \ref{4.7} vor und wir sind fertig. Es sei also $p < \rank \Theta'$.

Wir betrachten eine zusammenh\"angende offene $X$-Umgebung $X_{0}$ von $x_{0}$, auf der ein Erzeugendensystem $\theta^{(1)},\ldots,\theta^{(m)}$ von $\Theta'$ existiert. Wir d\"urfen annehmen:
\BegIT
   \item
      $\theta^{(1)}|_{x},\ldots,\theta^{(p)}|_{x}$ sind linear unabh\"angig f\"ur alle $x \in X_{0}$.
\EndIT
Es sei $\5D'_{p}$ die Menge aller Bl\"atter von $\5D'|_{X_{0}}$ der Dimension $p$. Dann ist
$$Z := \{x \in X_{0} : \dim \Theta'|_{x} \leq p\} = \{x \in X_{0} : \dim \Theta'|_{x} = p\}$$
eine analytische Teilmenge von $X_{0}$ und disjunkte Vereinigung der $C \in \5D'_{p}$. Deshalb induziert $\5D'_{p}$ eine quasi-analytische Struktur $\5T$ auf $Z$. Wir betrachten im Weiteren die Schichtung $(Z,\5T)$.

 Sei 
$\TT_{x} := \sum_{\nu=1}^p \C\theta^{(\nu)}|_{x}$ f\"ur $x \in X_{0}$. Dann ist
$\TT(Z,x)=\TT_{x}$ f\"ur alle $x \in Z$.

Im Weiteren bezeichnen wir mit $|\,.\,|$  die euklidische Norm und mit $B_{1}$ die offenen Kugeln im $\C^p_{z}$ und mit $B_{2}$ die offenen Kugeln im $\C^q_{w}$ (dabei sei $q := n-p$). Au{\ss}erdem machen wir gelegentlich Gebrauch von den \"ublichen Identifikationen
$$\TT'_{x} = \sum_{\nu=1}^p \C\frac{\partial}{\partial z_{\nu}}\bigg|_{x} = \C^p_{z}, \quad
\TT''_{x} = \sum_{\nu=1}^q \C\frac{\partial}{\partial w_{\nu}}\bigg|_{x} = \C^q_{w} \quad \Text{f\"ur} x = (z,w) \in \C^n.$$

Sei $B(x,\eps_{1},\eps_{2}) := B_{1}(z,\eps_{1}) \times B_{2}(w,\eps_{2})$ f\"ur $x = (z,w) \in \C^n$, 
$\eps_{1},\eps_{2} >0$.
Wir d\"urfen annehmen:
\BegIT
   \item
      $X_{0} \Opss \C^n$, $x_{0} = 0$ und $X_{0}$ ist von der Form $X_{0} = B(0,\eps_{1},\eps_{2})$,
   \item
      $A_{0} := B_{1}(0,\eps_{1}) \times \{0\}$ ist ein $\5D'$-Pl\"attchen von $A$
   \item f\"ur alle $x \in A_{0}$ ist
      $\theta^{(\nu)}\big|_{x} = \cases{\frac {\partial} {\partial z_{\nu}}|_{x} & $\nu = 1,\ldots,p$,\cr\cr
                                                               0 & $\nu = p+1,\ldots,m$.}$
   \item
      $\TT_{x} \cap \TT''_{x} = \{0\}$ f\"ur alle $x \in X_{0}$.
\EndIT

F\"ur die Projektion $\Mapping \pi {\C^n}{\C^p_{z}}$ gilt dann: $\Mapping \pi {\TT_{x}}{\TT'_{z}}$ ist --bei Benutzung der \"ublichen Identifikationen-- f\"ur alle $x \in X_{0}$ ein Isomorpismus und $\Mapping \pi Z {B_{1}(0,\eps_{1})}$ ist
bez\"uglich $\5T$ lokal biholomorph.

Wir betrachten den Isomorphismus $\Mapping \pi {\TT_{x}}{\TT'_{z}}$ genauer:

F\"ur $\nu = 1,\ldots,p$ ist
$$\theta^{(\nu)} = \frac{\partial}{\partial z_{\nu}} + \eta^{(\nu)}_{1}+ \eta^{(\nu)}_{2}$$
mit
$$\eta^{(\nu)}_{1}= \sum_{\lambda=1}^p a^{(\nu)}_{\lambda}\frac{\partial}{\partial z_{\lambda}}, \quad
\eta^{(\nu)}_{2}= \sum_{\mu=1}^q b^{(\nu)}_{\mu}\frac{\partial}{\partial w_{\mu}}.$$

Sei $\zeta = \sum_{\nu=1}^p \zeta_{\nu} \frac{\partial}{\partial z_{\nu}}$ und 
$\theta = \sum_{\nu=1}^p c_{\nu}\theta^{(\nu)}\big|_{x} \in \TT_{x}$ der Vektor mit $\pi(\theta)=\zeta$.
Dann ist
$$  \pmatrix{c_{1} \cr \vdots \cr c_{p}}  = N \cdot \pmatrix{\zeta_{1}\cr \vdots \cr \zeta_{p}}\quad$$
mit einer auf $X_0$ holomorphen Matrix $N$.

Weil $N(0)$ die Einheitsmatrix ist, d\"urfen wir annehmen:

\BegIT
   \item $\big|(c_{1},\ldots,c_{p})\big| \leq 2|\zeta|$.
\EndIT
Au{\ss}erdem d\"urfen wir noch annehmen:
\BegIT
   \item $\Big|\eta_{2}^{(\nu)}\big|_{x}\Big| < \frac 1p$ f\"ur $\nu = 1,\ldots, p$ und $x \in U_{0}$.
\EndIT
F\"ur die Komponente $\sum_{\nu=1}^p c_{\nu}\eta^{(\nu)}_{2}$ von $\theta$ gilt dann:
$$\Big|\sum_{\nu=1}^p c_{\nu} \eta^{(\nu)}_{2}\Big| \leq 2|\zeta| \quad  \Text{\"uberall auf} X_{0} \ \cap\  \{z\} \times \C^q_{w}.$$
Aus beweistechnischen Gr\"unden denken wir uns $\eps_{1},\eps_{2}$ so gew\"ahlt, da{\ss}
\BegIT
   \item $2\eps_{1} < \eps_{2}$
\EndIT
Die folgenden Betrachtungen zu dem lokalen Biholomorphismus $\Mapping \pi Z {B_{1}(0,\eps_{1})}$ sind f\"ur das Weitere wichtig.

Sei $\Bul x = (\Bul z, \Bul w) \in Z$. Dann existiert ein $\5T$-Pl\"attchen $\Bul C$ zu $\Bul x$ derart, da{\ss}
$\Mapping \pi {\Bul C}{\Bul U := \pi(\Bul C)}$ eine biholomorphe Abbildung ist. Wir d\"urfen annehmen, da{\ss} $\Bul U$ eine Kugel $B_{1}(\Bul z,\rho)$ ist. Wir nennen $\Mapping \pi {\Bul C}{B_{1}(\Bul z,\rho)}$ dann eine
\DEF{$\5T$-Kugel}\index{T-Kugel@$\5T$-Kugel} mit Mittelpunkt $\Bul x$.

Seien $\Mapping {\gamma, \tilde\gamma}{I}{Z}$ zwei $\5T$-Wege (dabei sei $I \subset \R$ ein Intervall) mit 
$\pi \circ \gamma = \pi \circ \tilde\gamma$. Gibt es einen Punkt 
$\Bul t \in I$ mit $\gamma(\Bul t) = \tilde \gamma(\Bul t)$, so ist $\gamma = \tilde \gamma$. Wir nennen $\gamma$ die \DEF{$\5T$-Liftung}\index{T-Liftung@$\5T$-Liftung} von $\pi \circ \gamma$ durch den Punkt 
$\Bul x = \gamma(\Bul t)$.

Doch nun weiter! Wir bezeichnen im Folgenden $\C^p_{z}$-Komponenten mit dem Index 1 und 
$\C^q_{w}$-Komponenten mit dem Index 2.

Wir betrachten ein konstantes Vektorfeld
$$\zeta = \sum_{\nu = 1}^p \zeta_{\nu}\frac \partial {\partial z_{\nu}},\ (\zeta_{1},\ldots,\zeta_{p}) \in \C^p,
 |\zeta| = 1$$
auf $B_{1}(0,\eps_1)$. Zu $\zeta$ existiert ein holomorphes Vektorfeld $\theta$ auf $X_{0}$ mit 
$\pi(\theta|_{x}) = \zeta|_{z}\ \forall x = (z,w) \in X_{0}$. Es gilt: $\theta_{1} = \zeta$, insbesondere 
$\big|\theta_{1}|_{x}\big| = 1$, und $\big|\theta_{2}|_{x}\big| \leq 2$ f\"ur alle $x \in X_{0}$. 

Wir w\"ahlen jetzt $\delta_{1},\delta_{2}$ mit $0 < \delta_{1} < \eps_{1}/4$ und $\delta_{1}< \delta_{2} < \eps_{2}/4$.
Wir betrachten einen Punkt $\Bul x = (\Bul z,\Bul w) \in Z \cap B(0,\delta_{1},\delta_{2})$ und dazu ein konstantes Vektorfeld $\zeta$ sowie das zugeh\"orige Vektorfeld $\theta$ wie oben. Sei $\Mapping \gamma I {X_{0}}$ die maximale L\"osung der (reellen) Differentialgleichung $\gamma' = \theta \circ \gamma$ mit $\gamma(0) = \Bul x$.
$I$ ist ein offenes Intervall in $\R$. Mit Standard{\"u}berlegungen zeigt man, da{\ss} $[0,\delta_{1}[\  \subset I$.

Wir betrachten eine $\5T$-Kugel $\Mapping \pi {\Bul C}{B_{1}(\Dot z,\rho)}$ mit $0 < \rho \leq \delta_{1}$.
F\"ur den Weg
$$\Mapping{\tilde \gamma}{[0,\rho[\ }{\Bul C}, \quad \tilde\gamma(t) := 
(\pi|_{\Dot C})^{-1}\big(\Dot z + t(\zeta_{1},\ldots,\zeta_{p})\big)$$
gilt dann $\pi\big(\tilde\gamma'(t)\big) = \zeta$, also $\tilde \gamma'(t) = \theta \circ \tilde\gamma(t)$ und
 $\tilde\gamma(0) = \Bul x$. Folglich ist $\tilde\gamma = \gamma|_{[0,\rho[\ }$. Der Weg $\gamma|_{[0,\rho[\ }$ 
 verl\"auft in $\Dot C$. Weil $Z$ $X_{0}$-abgeschlossen ist, folgt $\gamma(\rho) \in Z$. 
 Ist nun $t_{*} \in ]0,\delta_{1}]$ mit $x_{*} := \gamma(t_{*}) \in Z$, so folgt mit der selben Argumentation wie im Fall 
 $t = 0$, da{\ss} $\gamma$ auf einer $]0,\delta_{1}]$-offenen Umgebung von $t_{*}$ in einem $\5T$-Pl\"attchen 
 verl\"auft und da{\ss} die Werte von $\gamma$ in den Endpunkten in $Z$ liegen. Nun schlie{\ss}t man sofort, da{\ss} $\gamma$ ein in $Z$ verlaufender $\5T$-Weg ist; $\gamma$ ist die $\5T$-Liftung der Strecke
 $[\Dot z,\Dot z + \delta_{1}(\zeta_{1},\ldots,\zeta_{p})]$ durch~$\Dot x$. 
 Wir nennen $\Mapping \gamma {[0,\delta_{1}]}Z$ die \DEF{$\5T$-Strecke}\index{T-Strecke@$\5T$-Strecke}
 mit ``Startpunkt'' $\Dot x$,  ``Richtung'' $\zeta$ und ``L\"ange''~$\delta_{1}$.
 
 Wir halten fest: Ist $\Mapping \gamma {[0,\delta_{1}]}Z$ eine $\5T$-Strecke wie oben und trifft $\gamma$ eine $\5T$-Kugel $\Mapping \pi {C_{*}}{B_{1}(z_{*},\rho)}$, so gilt
$$\gamma(t) = (\pi|_{C_{*}})^{-1}(\Dot z + t\zeta)\quad \forall t \in \R \Text{mit} \Dot z + t\zeta \in B_{1}(z_{*},\rho).$$
Nun sei
$$\Dot C := \big\{x \in Z  : \Text{es gibt eine in $\Dot x$ startende $\5T$-Strecke 
$\Mapping \gamma {[0,\delta_{1}]}Z$
und $t \in [0,\delta_{1}[\ $ mit $\gamma(t) = x$}\big\}.$$
Man zeigt nun mit leichter M{\"u}he:
\BegEN
   \item \label{BEH1}$\Dot C$ ist $\5T$-offen,
   \item \label{BEH2}$\Dot C$ ist eine analytische Teilmenge von $\Dot U := B\big((\Dot z,0),\delta_{1},\eps_{2}\big)$,
   \item \label{BEH3}$\Dot C$ ist die $\5D'_{p}$-Komponente von $Z \cap \Dot U$, welche $\Dot x$ enth\"alt.
\EndEN
\medskip

Nun kommen wir zum entscheidenden Argument f\"ur den Beweis unseres Satzes \ref{4.9}. 

Sei 
$\Dot x \in Z \cap B(0,\delta_{1}/2,\delta_{2})$. 
Dann ist $B(0, \delta_{1}/2,\eps_{2}) \subset B(\Dot z,\delta_{1}, \eps_{2})$. Also sind alle $\5D'$-Komponenten von
$A \cap B(0,\delta_{1}/2,\delta_{2})$ analytische Mengen in $B(0,\delta_{1}/2,\delta_{2})$, insbesondere bez\"uglich der Relativtopologie abgeschlossen. \Qed

Im Fall einer regul\"aren Bl\"atterung $\5F$ ist die Spallek-Zerlegung $\5D'_{\5F}$ von $\Theta_{\5F}$ nat\"urlich gleich der Zerlegung~$\5D_{\5F}$.

Sei $\5F$ eine $p$-dimensionale (koh\"arente) holomorphe Bl\"atterung
auf $X$; Singularit\"aten sind zugelassen. $\5F$ ist definiert durch eine vollst\"andige involutive koh\"arente $\5O$-Untergarbe $\Theta_{\5F}$ von $\Theta$ vom Rang $p$ bzw.~eine vollst\"andige involutive koh\"arente $\5O$-Untergarbe $\Omega_{\5F}$ von $\Omega$ vom Rang $q:=n-p$. Die Garben $\Theta_{\5F}$ und $\Omega_{\5F}$ sind dual zueinander. Es ist
$\Sing \5F  = \Sing \Theta_{\5F} = \Sing\Omega_{\5F}$ eine analytische Menge der Kodimension $\geq 2$.
Die Einschr\"ankung $\5F|_{(X \sm \Sing\5F)}$ ist eine regul\"are holomorphe Bl\"atterung.
Die Spallek-Zerlegung $\5D'_{\5F}$ von $\Theta_{\5F}$ nennen wir auch die
\DEF{Spallek-Zerlegung von $\5F$}\index{Spallek-Zerlegung!von $\5F$}.

Neben $\5D'_{\5F}$ kann in geeigneten F\"allen eine weitere mit $\5F$ zusammenh\"angende quasi-analytische Zerlegung von $X$ eingef\"uhrt werden, vgl.~\cite{Reiffen1} und \cite{Reiffen2}.
Dazu fundamental ist der Begriff der Stammfunktion: Sei $U \Opss X$. Eine Funktion $f \in \5O(U)$ hei{\ss}t
\DEF{$\5F$-Stammfunktion}\index{F-Stammfunktion@$\5F$-Stammfunktion}, wenn $df \in \Omega_{\5F}(U)$ gilt. Die
$\5F$-Stammfunktionen bilden die Untergarbe $\5O_{\5F}$ von $\5O$.

Eine Schicht $(Z,\5T)$ von $X$ hei{\ss}t eine 
\DEF{$\5F$-Integralvariet\"at}\index{F-Integralvariet\"at@$\5F$-Integralvariet\"at}%
\index{Integralvariet{\"a}t}, wenn f\"ur alle $z \in Z$ und $f \in (\5O_{\5F})_{z}$ die Beschr\"ankung des Keimes $f$ auf 
$(Z,\5T)$ konstant ist.

\Beg{Satz}\label{4.10}
   Jedes Blatt $A \in \5D'_{\5F}$ ist eine $\5F$-Integralvariet\"at.
\End

{\bf Beweis:} Sei $A_{0}$ ein Pl\"attchen von $A$, $x \in A_{0}$, $f \in (\5O_{\5F})_{x}$. Dann gilt:
$df \in (\Omega_{\5F})_{x}$, $df \big((\Theta_{\5F})_{x}\big) = 0$, $d(f|_{A_{0}}) = 0$, $f|_{A_{0}} = \const$. \Qed

\Beg{Definition}\label{4.11}
   Eine zusammenh\"angende lokal-analytische Teilmenge $A$ von $X$ hei{\ss}t ein
   \DEF{lokales $\5F$-Blatt}\index{lokales F-Blatt@lokales $\5F$-Blatt}, wenn gilt:
   \BegEN
      \item $A$ ist eine $\5F$-Integralvariet\"at\label{4.11.1}
      \item $\dimMin A \geq p$ $(=\dim \5F)$\label{4.11.2}
      \item Ist $B$ eine weitere lokal-analytische Teilmenge von $X$ mit den Eigenschaften \ref{4.11.1} und
      \ref{4.11.2}, so gilt $B_{x} \subset A_{x}$ f\"ur alle $x \in A \cap B$.\label{4.11.3}
   \EndEN
\End

Sei $\5A_{\5F}$ die Menge der lokalen $\5F$-Bl\"atter und $$X_{\5F} := \bigcup_{A \in \5A_{\5F}} A.$$
$\5A_{\5F}$ erf\"ullt die Voraussetzungen von \ref{1.10}, definiert deshalb eine quasi-analytische Struktur 
auf $X_{\5F}$. Die Schichten von $(X_{\5F},\5T)$ heissen \DEF{$\5F$-Bl\"atter}\index{Blatt}. 
Sie sind $\5F$-Integralvariet\"aten. Im Allgemeinen ist $X_{\5F} \not= X$. Wenn $X_{\5F} = X$, so sagen wir
\DEF{$\5F$ hat Bl\"atter \"uberall}\index{Bl\"atter \"uberall} und bezeichnen mit 
$\5D_{\5F}$\index{DF@$\5D_{\5F}$} die Menge der $\5F$-Bl\"atter \DEF{(Bl\"atterraum)}\index{Bl\"atterraum}.
 Ist $\5F$ regul\"ar, so stimmt dieser 
Bl\"atterraum mit dem in \ref{4.7} genannten Bl\"atterraum $\5D_{\5F}$ \"uberein. 
Im Allgemeinen ist $\5D_{\5F}' \not= \5D_{\5F}$, auch wenn $\5F$ Bl\"atter \"uberall hat:

\Beg{Beispiel}\label{4.12}
   Sei $\Mapping f {\C^2}\C$, $(z,w) \mapsto z \cdot w$, $\Omega' := \5O\,df = \5O(w\,dz + z\,dw)$. 
   Dann ist $\Omega'$ offenbar involutiv und  vollst\"andig, weil $\Sing \Omega' = \{0\}$ 
   2-codimensional ist. 
   Die Garbe $\Theta' = \5O\big(z\frac \partial{\partial z} - w \frac \partial{\partial w}\big)$ ist die Annullatorgarbe
   von $\Omega'$. Die Garben $\Omega'$ bzw.~$\Theta'$ definieren eine 1-dimensionale holomorphe Bl\"atterung
   $\5F$ auf $X$. Es besteht $\5D_{\5F}'$ aus den Bl\"attern $\{0\}$, $\C_{z}^* \times \{0\}$, $\{0\} \times \C_{w}^*$
   sowie $\{(z,w) \in \C^2 : zw = c\}$, $c \in \C^*$.
   $\5F$ hat Bl\"atter \"uberall und $\5D_{\5F}$ besteht aus den Bl\"attern $\{(z,w) \in \C^2 : zw = c\}$, $c \in \C$.
\End

{\bf Beweis:} Sei $x_{0} := (z_{0},w_{0})$. Dann ist $\Theta'|_{x_{0}}=\C \cdot (z_{0},-w_{0})$, woraus sofort die Aussage \"uber $\5D_{\5F}'$ folgt. Auf $\C^2 \sm \{0\}$ ist $\5F$ regul\"ar und besitzt dieselben Bl\"atter wie 
$\5D_{\5F}'|_{\C^2 \sm \{0\}}$, also die Bl\"atter $\{(z,w) \in \C^2 : z \cdot w = c\}$, $c \in \C^*$, sowie
$\C_{z}^* \times \{0\}$ und $\{0\} \times \C^*_{w}$. Nehmen wir nun $0$ hinzu, so ist offenbar die
 $X$-abgeschlossene H\"ulle von $\big(\C^*_{z} \times \{0\}\big) \cup \big(\{0\} \times \C^*_{w}\big)$, das ist
 $\big(\C_{z} \times \{0\}\big) \cup \big(\{0\} \times \C_{w}\big)$, eine $\5F$-Integralvariet\"at. Weil \ref{4.11} erf\"ullt ist, handelt es sich um ein lokales $\5F$-Blatt und dann sofort um ein $\5F$-Blatt. \Qed

\Beg{Satz}\label{4.13}
   Sei $\5F$ eine $p$-dimensionale holomorphe Bl\"atterung mit Bl\"attern \"uberall. Dann ist die 
   quasi-analytische Zerlegung $\5D_{\5F}$ von $X$ schwach-analytisch von der Feinheit $\leq p$.
\End

{\bf Beweis:} Sei $A \in \5D_{\5F}$ und $x_{0} \in A$. Wir betrachten ein lokales $\5F$-Blatt $A_{0}$ mit 
$x_{0} \in A_{0} \subset A$ und zeigen $A_{0} \in \5A_{p}(A)$. Sei $B \in \5I_{p}(A)$, $B \subset A$.
Wir zeigen, da{\ss} jede $X$-Zusammenhangskomponente $B'$ von $B$ eine $\5F$-Integralvariet\"at ist. Dann folgt mit \ref{4.11}.\ref{4.11.3}, da{\ss} $B_{z} \subset (A_{0})_{z}$ f\"ur alle $z \in A_{0} \cap B$ und wir sind fertig. 

Wir d\"urfen annehmen, da{\ss} $B$ zusammenh\"angend ist. Sei $x \in B$ und $f \in (\5O_{\5F})_{x}$. Wir betrachten eine zusammenh\"angende offene $X$-Umgebung $U$ von $x$ mit den folgenden Eigenschaften:
$B \cap U$ ist zusammenh\"angend und $f \in \5O_{\5F}(U)$. Ist $C$ eine $\5D_{\5F}$-Komponente von 
$A \cap U$, so ist $f$ auf $C$ und dann auch auf $B \cap U \cap C$ konstant. Da $A \cap U$ h\"ochstens 
abz\"ahlbar viele $\DF$-Komponenten besitzt, ist $N := f(B \cap U)$ eine h\"ochstens abz\"ahlbare Teilmenge 
von $\C$.
Da $N$ zugleich auch zusammenh\"angend ist, mu{\ss} $N$ einelementig sein, das hei{\ss}t $f|_{B \cap U}$ ist konstant. Also ist $B$ eine Integralvariet\"at. \Qed

\Beg{Definition}
   Sei $\5F$ eine $p$-dimensionale holomorphe Bl\"atterung. Ein $\5F$-Blatt hei{\ss}t \DEFI{stark}, wenn f\"ur jedes lokale $\5F$-Blatt $A_{0} \subset A$ gilt:
   \begin{quote}
      Ist $B$ eine lokal-analytische Teilmenge von $X$ und eine $\5F$-Integralvariet\"at, so ist 
      $B_{x} \subset (A_{0})_{x}$ f\"ur alle $x \in A_{0} \cap B$ (beachte, da{\ss} hier im Vergleich zu 
      \ref{4.11}.\ref{4.11.2} die Forderung Ò$\dimMin B \geq p$Ó weggelassen wurde).
   \end{quote}
   Die Bl\"atterung $\5F$ hei{\ss}t  \DEFI{stark}, wenn $\5F$ \"uberall starke Bl\"atter besitzt.
\End

\Beg{Satz}\label{4.15}
   Sei $\5F$ eine starke $p$-dimensionale Bl\"atterung. Dann ist $\DF$ schwach-analytisch von der Feinheit $0$.
   Zu jedem $A \in \DF'$ gibt es ein $B \in \DF$ mit $A \subset B$.
\End

{\bf Beweis:} Der erste Teil der Aussage folgt mit einer leichten Modifikation des Beweises von \ref{4.13}. Sei 
$A \in \DF'$. Wir betrachten ein $\DF'$-Pl\"attchen $\0 A$ von $A$. Sei $\0 x \in \0 A$ und $B \in \DF$ mit 
$\0 x \in B$. Sei $\0B$ ein lokales $\5F$-Blatt mit $\0x \in \0B \subset B$. Wegen \ref{4.10} d\"urfen wir
$\0A \subset \0B$ annehmen. Also mu{\ss} $A \cap B$ eine $\DF'$-offene Teilmenge von $A$ sein. Wir nehmen an da{\ss} $A \cap B \neq A$. Sei $\1x \in A$ ein $\DF'$-Randpunkt von $A \cap B$. Wir betrachten ein 
$\DF'$-Pl\"attchen $\1A \subset A$ mit $\1x \in \1A$ sowie ein $C \in \DF$ mit $\1x \in C$. Mit obigem Argument gibt es ein lokales $\5F$-Blatt $\1C$ mit $\1A \subset \1C \subset C$. In $\1A$ liegen Punkte aus $A \cap B$. Es folgt:
$B \cap C \neq \emptyset$, $B = C$, $\1x \in B$. Widerspruch! Also mu{\ss} $A \subset B$ gelten. \Qed

Wir k\"onnen die zweite Aussage in \ref{4.15} auch so formulieren:

\Beg{Satz}\label{4.16}
   Sei $\5F$ eine starke $p$-dimensionale Bl\"atterung. Dann definiert $\DF'$ auf jedem $B \in \DF$ eine Schichtung.
\End

Ist $\5F$ eine 1-dimensionale Bl\"atterung mit Bl\"attern \"uberall, so ist $\5F$ offenbar automatisch stark. Das gilt insbesondere f\"ur unser Beispiel \ref{4.12}.

F\"ur quasi-analytische Zerlegungen von $X$ sind zwei Begriffe von Stammfunktion sinnvoll:

\Beg{Definition}
   Sei $\5D$ eine quasi-analytische Zerlegung von $X$, $U \Opss X$ und $f \in \5O(U)$. Dann hei{\ss}t $f$ eine
   \DEF{$\5D$-Stammfunktion}\index{D-Stammfunktion@$\5D$-Stammfunktion}\index{Stammfunktion}, wenn
   f\"ur jede $X$-Zusammenhangskomponente $U'$ von $U$ die Einschr\"ankung $f|_{U'}$ auf den Bl\"attern
   von $\5D|_{U'}$ konstant ist, oder anders formuliert, wenn $f$ bez\"uglich der $\5D$-Topologie lokal konstant ist. 
   $f$ hei{\ss}t eine \DEFI{differentielle $\5D$-Stammfunktion}\index{Stammfunktion!differentielle}, wenn
   $df \big ((\Theta^{\5D})_{x}\big) = 0$ f\"ur alle $x \in U$.
   Die $\5D$-Stammfunktionen bilden die Garbe $\5O^{\5D}$, die differentiellen $\5D$-Stammfunktionen die 
   Garbe $\5O^{\5D}_{d}$\index{OD@$\5O^{\5D}$}\index{ODd@$\5O^{\5D}_{d}$}.
\End

\Beg{}\label{4.18}
   Sei $\5D$ eine quasi-analytische Zerlegung von $X$. Dann ist $\5O^{\5D} \subset \5O^{\5D}_{d}$. 
\End

{\bf Beweis:} Sei $U \Opss X$, $f \in \5O(U)$ eine $\5D$-Stammfunktion und $x \in U$. Sei $B$ ein $\5D$-Pl\"attchen mit $x \in B \subset U$. Wir betrachten $g := f - f(x)$. Dann ist $g$ eine auf $B$ verschwindende holomorphe Funktion. Deshalb gilt: $dg = df$ annulliert alle Tangentialvektorr\"aume $\TT(B,y)$, $y \in B$. Daraus folgt unsere Aussage. \Qed

Da allgemeine quasi-analytische Zerlegungen sehr ÒchaotischÓ sein k\"onnen, ist in \ref{4.18} die Gleichheit im Allgemeinen nicht gegeben:

\Beg{Beispiel}
   Sei $X = \C^2 = \C_{z} \times \C_{w}$. F\"ur $x = (z,w) \in X$ sei das Blatt $\5D(x)$ wie folgt definiert:
   $$\5D(x) := \cases{\{x\} & wenn $z$ rational ist, d.h.~$z \in \7Q^2$, \cr
                                  \{z\} \times \C_{w} & sonst.}$$
   Das System der $\5D(x)$, $x \in X$, definiert eine quasi-analytische Zerlegung $\5D$ auf $X$. Dabei ist 
   $\Theta^{\5D}$ die Nullgarbe.                               
Sei $\1D$ ein Gebiet in $\IND \C z$, $\2D$ ein Gebiet in $\IND \C w$ und $D := \1D \times \2D$. Dann ist
$$\5O^{\5D}(D) = \{f \in \5O(D) : \frac {\partial f}{\partial w} = 0\} = \5O(\1D), \quad \5O^{\5D}_{d}(D) = \5O(D).$$
\End

\Beg{}
   Sei $\Theta'$ eine involutive koh\"arente $\5O$-Untergarbe von $\Theta$ und $\5D'$ die zu $\Theta'$ 
   geh\"orige Spallek-Zerlegung. Dann ist $\5O^{\5D'} = \5O^{\5D'}_{d}$.
\End

{\bf Beweis:} Zun\"achst ist $\Theta' \subset \Theta^{\5D'}$. Nun sei $U \Opss X$ und $f \in \5O^{\5D'}_d(U)$.
Dann gilt insbesondere: $df(\Theta'_{x})=0$ f\"ur alle $x \in U$. Dann mu{\ss} $f$ auf den in $U$ liegenden
$\5D'$-Pl\"attchen konstant sein. \Qed

\Beg{}\label{4.21}
   Sei $\5F$ eine $p$-dimensionale holomorphe Bl\"atterung mit Bl\"attern \"uberall und $\5D := \DF$, $\5D':=\DF'$.
   Dann gilt:
   $$\5O^{\5D} = \5O_{\5F} = \5O^{\5D}_{d} = \5O^{\5D'} = \5O^{\5D'}_{d}.$$
\End
{\bf Beweis:} Weil jedes $\5F$-Blatt eine $\5F$-Integralvariet\"at ist, gilt: $\5O_{\5F} \subset \5O^{\5D}$.
Sei $U \Opss X$ und $f \in \5O^{\5D}_{d}(U)$. Auf $U^* := U \sm \Sing \5F$ gilt dann (weil $\Omega_{\5F}$ und
$\Theta_{\5F}$ zueinander dual sind): $df|_{U^*} \in \Omega_{\5F}(U^*)$. Weil $\Omega_{\5F}$ vollst\"andig ist, folgt: $df \in \Omega_{\5F}(U)$, also $f \in \5O_{\5F}(U)$. Also gilt auch $\5O^{\5D}_{d} \subset \5O_{\5F}$.
Weil jedes $A \in \5D'$ eine $\5F$-Integralvariet\"at ist, gilt: $\5O_{\5F} \subset \5O^{\5D'}$. Sei $U \Opss X$ und
$f \in \5O^{\5D'}_{d}(U)$. Dann schlie{\ss}t man wie oben, da{\ss} $f \in \5O_{\5F}(U)$. Also gilt auch
$\5O^{\5D'}_{d} \subset \5O_{\5F}$. \Qed

\subsection{Durch Abbildungen definierte quasi-analytische Zerlegungen und holomorphe Bl\"atterungen%
}\label{AdZ}

Wir beginnen mit einigen allgemeinen Betrachtungen, die an unser Beispiel \ref{4.6} ankn\"upfen. 

Sei $\Mapping f X W$ eine holomorphe Abbildung in den komplexen Raum $W$. Dann bezeichnen wir mit $\Df$
die durch $f$ auf $X$ definierte quasi-analytische Zerlegung, also das System der Niveaumengen von $f$. 
Sei $\Omega'_{f} := f^*(\Omega_{W})$. Ist $W$ eine lokal-analytische Teilmenge eines $\C^N$ (bei einer lokalen Betrachtung kann man davon stets ausgehen) und $f = (f_{1},\ldots,f_{N})$, so ist 
$\Omega'_{f} = \sum_{j=1}^N \5O\,df_{j}$.
Wir nennen $\Sing f := \Sing \Omega_{f}$ die \DEFI{Singularit\"atenmenge} von $f$. Sie ist eine niederdimensionale analytische Teilmenge von $X$. 
Sei $\Omega_{f} := \widetilde{\Omega'_{f}}$ die Vervollst\"andigung von $\Omega'_{f}$ und 
$\Theta_{f} := \Omega^0_{f}$ die Annullatorgarbe von $\Omega_{f}$. Dann definieren $\Omega_{f}$ und 
$\Theta_{f}$ dieselbe holomorphe Bl\"atterung $\5F_{f}$ auf $X$; diese hat die Dimension $p=n-q$, dabei ist
$q = \Rang \Omega_{f} = \Rang \Omega_{f}'$, $p = \Rang \Theta_{f}$. In \cite{Reiffen1} wird $f$ ein schwaches Integral von $\5F_{f}$ genannt.

Sei $\0x \in X \sm \Sing f$. Wir betrachten eine zusammenh\"angende offene $X$-Umgebung $U$ von $\0x$, 
$U \subset X \sm \Sing f$, sowie $\Mapping f U W$. Wir d\"urfen annehmen, da{\ss} $W$ eine lokal-analytische Teilmenge eines $\C^N$ ist. Dann hat die Abbildung $\Mapping{f=(f_{1},\ldots,f_{N})}U{\C^N}$ \"uberall den konstanten Rang $q$. Aufgrund des Rangtheorems (vgl. etwa \cite[Seite 313]{ReiffenTrapp} oder 
\cite[Seite 18]{Narasimhan}) d\"urfen wir sogar von der Situation \ref{4.5} ausgehen, 
d.h.~$U$ ist ein Gebiet im $\C^n$ von der Form
$U = \1D \times \2D$, wobei $\1D \subset \C^p_{z}$, $\2D \subset \C^q_{w}$ Gebiete sind, und es ist $W = \2D$,
$f(z,w) = w$. 

\Beg{Satz}\label{5.1}
   In der oben beschriebenen Situation gilt auf $X \sm \Sing f$: $\5F_{f}$ ist eine regul\"are holomorphe
   Bl\"atterung, deren zugeh\"orige quasi-analytische Zerlegung dort mit
   $\Df$ \"ubereinstimmt.
\End

Aus \ref{5.1} und der Vollst\"andigkeit der Garbe $\Omega_{f}$ bzw.~$\Theta_{f}$ folgt:

\Beg{}\label{5.2}
  Ist $\Mapping {\Dot f} X {\Dot W}$ eine weitere holomorphe Abbildung und gilt $\Df = \5D_{\Dot f}$, so ist
  $\5F_{f} = \5F_{\Dot f}$.
\End

\Beg{Definition}\label{5.3}
   Seien $\5D$ eine quasi-analytische Zerlegung
   und $\5F$ eine holomorphe Bl\"atterung von $X$. 
   Unter einer \DEF{lokalen Beschreibung von $\5D$ bzw. $\5F$}%
   \index{lokale Beschreibung} verstehen wir eine holomorphe Abbildung $\Mapping f U W$ einer
   zusammenh\"angenden $X$-offenen Teilmenge von $X$ in einen komplexen Raum $W$ derart, da{\ss}
   $\5D|_{U} = \5D_{f}$ bzw.~$\5F|_{U} = \5F_{f}$ ist. $\5D$ bzw.~$\5F$  hei{\ss}t \DEFI{abbildungsdefiniert}, 
   wenn es zu jedem Punkt $x \in X$ eine 
   lokale Beschreibung $\Mapping f U W$ von $\5D$ bzw.~$\5F$ gibt mit $x \in U$.
\End

\Beg{Satz}\label{5.5}
   Sei $\5D$ eine abbildungsdefinierte quasi-analytische Zerlegung von $X$. Dann definiert $\5D$ in nat\"urli\-cher
   Weise eine abbildungsdefinierte holomorphe Bl\"atterung auf $X$.
\End

{\bf Beweis:} Seien $\Mapping f U W$ und $\Mapping {\Dot f} {\Dot U} {\Dot W}$ lokale Beschreibungen von $\5D$.
Wegen \ref{5.2} gilt auf $U \cap \Dot U$: $\5F_{f}|_{U \cap \Dot U} = \5F_{\Dot f}|_{U \cap \Dot U}$.
Deshalb k\"onnen die lokal definierten Bl\"atterungen $\5F_{f}$ zu einer holomorphen Bl\"atterung $\5F$ auf $X$ verklebt werden. $\5F$ ist nach Konstruktion abbildungsdefiniert. \Qed

\Beg{Definition}\label{5.6}
    Eine quasi-analytische Zerlegung $\5D$ von $X$ hei{\ss}t im Punkte $x \in X$ \DEFI{regul\"ar}, wenn es
    eine lokale Beschreibung $\Mapping f U {\C^q}$ von $\5D$  mit $x \in U$ gibt, bei der $f$ eine holomorphe
    Submersion ist. Die $X$-offene Menge
    $$\Reg \5D := \{x \in X : \5D \Text{ ist in} x \Text{regul\"ar}\kern-3pt\}$$
    hei{\ss}t die \DEFI{regul\"are Menge} von $\5D$ und $\Sing \5D  := X \sm \Reg \5D$ die
    \DEFI{Singularit\"atenmenge} von $\5D$. Wenn $\Reg \5D = X$ ist,  hei{\ss}t $\5D$ \DEFI{regul\"ar}.
\End

Die regul\"aren Zerlegungen entsprechen offensichtlich genau den regul\"aren Bl\"atterungen.

Eine Teilmenge $M$ von $X$ hei{\ss}t bekanntlich \DEFI{analytisch d\"unn}, wenn es zu jedem Punkt $x \in M$ 
eine offene $X$-Umgebung $U$ von $x$ und eine niederdimensionale analytische Menge $A$ in $U$ gibt mit
$M \cap U \subset A$. Aus \ref{5.1} folgt:

\Beg{Satz}\label{5.7}
   Die quasi-analytische Zerlegung $\5D$ von $X$ sei abbildungsdefiniert. Dann ist $\Sing \5D$ analytisch d\"unn.
\End

\Beg{Satz}\label{5.8}
   Die quasi-analytische Zerlegung $\5D$ von $X$ sei rein $p$-dimensional und $\Theta^{\5D}$ sei regul\"ar vom Rang $p$. Dann ist $\5D$ regul\"ar.
\End

Zum {\bf Beweis} vgl.~\cite[Satz 6]{Knoche}

\Beg{Definition}\label{5.9}
   Eine holomorphe Abbildung $\Mapping f X W$ von $X$ in einen komplexen Raum $W$ hei{\ss}t
   \DEFI{fasertreu}, wenn f\"ur alle Gebiete $U \Opss X$ und alle $g \in \5O(U)$ gilt: 
   \\
   ist $g|_{U \sm \Sing f}$ eine Stammfunktion von $\5D_{f}|_{U \sm \Sing f}$, so ist  $g$ bereits eine Stammfunktion
   von $\5D_{f}|_{U}$.
\End

Wir wollen die Eigenschaft, da{\ss} $g|_{U \sm \Sing f}$ eine Stammfunktion von $\5D_{f}|_{U \sm \Sing f}$ ist,
n\"aher untersuchen. Dazu sei $\Theta'_f := \Theta^{\5D_{f}}$.

Es ist $\Theta'_{f}|_{X \sm \Sing f} = \Theta_{f}|_{X \sm \Sing f}$ und dann $\Theta'_{f} \subset  \Theta_{f}$, weil
$\Theta_{f}$ vollst\"andig ist.

Sei $\omega \in \Omega(U)$. Dann gilt
$$
  \omega\big((\Theta'_{f})_{x}\big) = 0 \ \forall x \in U \ \iff\ 
                        \omega\big((\Theta_{f})_{x}\big) = 0 \ \forall x \in U
                        \ \iff \  \omega \in \Omega_{f}(U).
                        $$

F\"ur $g \in \5O(U)$ gilt:
\BegEA
   g|_{U \sm \Sing f} \in \5O^{\5D_{f}}(U \sm \Sing f) &\iff&
     g|_{U \sm \Sing f} \in \5O^{\5D_{f}}_{d}(U \sm \Sing f) \Text{\quad (vgl. \ref{5.1})} \\
      &\iff&dg \in \Omega_{f}(U \sm \Sing f)\\
        &\iff& dg \in \Omega_{f}(U) \Text{\quad (weil $\Omega_{f}$ vollst\"andig ist)}\\
        &\iff& g \in \5O_{\5F_{f}}(U)\\
        &\iff& g \in \5O^{\5D_{f}}_{d}(U).
\EndEA

Damit haben wir bewiesen:

\Beg{Satz}\label{5.10}
   F\"ur die holomorphe Abbildung $\Mapping f X W$ sind folgende Aussagen \"aquivalent:
   \BegEN
      \item $f$ ist fasertreu,
      \item $\5O^{\5D_{f}} = \5O^{\5D_{f}}_{d}$,
      \item $\5O^{\5D_{f}} = \5O_{\5F_{f}}$.
   \EndEN
\End

In \cite{Reiffen1} findet man Beispiele f\"ur fasertreue Abbildungen.

\Beg{Satz}\label{5.11}
   F\"ur die holomorphe Abbildung $\Mapping f X W$ sei $\5D_{f}$ reindimensional. Dann ist $f$ fasertreu.
\End

Zum {\bf Beweis} vgl.~\cite[5.4]{Reiffen1}. Jede offene holomorphe Abbildung $\Mapping f X W$ erf\"ullt die Vorraussetzungen von \ref{5.11} und ist daher fasertreu (vgl.~\ref{4.6} ff).

Die holomorphe Abbildung $\Mapping f X W$ hei{\ss}t 
\DEF{\"uberall faktorisierend}\index{uberall faktorisierend@\"uberall faktorisierend}\index{faktorisierend}, wenn gilt:
f\"ur alle $x \in X$ und alle $g \in (\5O_{\5F_{f}})_{x}$ gibt es ein $h \in (\5O_{W})_{f(x)}$ mit $g = h \circ f$.

\Beg{Satz}\label{5.12}
   Die holomorphe Abbildung $\Mapping f X W$ sei \"uberall faktorisierend, Dann ist $f$ fasertreu.
\End

{\bf Beweis:} Es ist immer $\5O^{\5D_{f}} \subset \5O_{\5F_{f}}$. Sei $\0x \in X$, $g \in (\5O_{\5F_{f}})_{\0x}$ und dazu $h \in (\5O_{W})_{f(\0x)}$ mit $g = h \circ f$. Wir betrachten Repr\"asentanten $\Mapping g U \C$ und
$\Mapping h V \C$ der Keime $g,h$ auf zusammenh\"angenden offenen Umgebungen $U$ bzw.~$V$ von
$\0x$ bzw.~$f(\0x)$ mit $f(U) \subset V$ und $g = h \circ f$ auf $U$. Dann gilt offensichtlich
$g \in \5O^{\5D_{f}}(U)$. \Qed

Der folgende Satz stammt von {\sc Malgrange} (vgl.~5.21 und 5.22 in \cite{Reiffen1}):

\Beg{Satz}\label{5.13}
   F\"ur die holomorphe Abbildung $\Mapping f X W$ gelte: $W$ ist eine Mannigfaltigkeit der Dimension
   $q=\Rang \Omega_{f}$ und $\codim \Sing f \geq 2$. Dann ist $f$ \"uberall faktorisierend, insbesondere fasertreu.
\End

Da{\ss} nicht jede holomorphe Abbildung $\Mapping f X W$ fasertreu ist, zeigt das folgende Beispiel:

\Beg{Beispiel}\label{5.14}
   Sei $X = \C^3_{z,w,u}$, $W = \C^2_{s,t}$ und $\Mapping f X W$ definiert durch \\
   $(z,w,u) \mapsto (z,z \cdot w \cdot u)$. Die Funktionalmatrix $Df$ von $f$ ist
   $$Df = \pmatrix{1 & 0 & 0 \cr w \cdot u & z \cdot u & z \cdot w}$$
   also ist $\Sing f = \big(\{0\} \times \C^2\big) \cup \big(\C \times \{(0,0)\}\big)$.
   F\"ur $(s,t) \in \C^2$ gilt:
   $$f^{-1}(s,t) = \cases{\{s\} \times \big\{(w,u) \in \C^2 : w \cdot u =\frac ts\big\} & falls $ s\neq 0$, \cr \cr
                                      \{0\} \times \C^2 & falls $s=0$.}
   $$
   Die Funktion $\Mapping g {\C^3}\C$, $g(z,w,u) := w \cdot u$, ist auf $\C^3 \sm \big(\{0\} \times \C^2\big)$
   eine $\5D_{f}$-Stammfunktion, nicht aber auf $\C^3$.
\End

Gem\"a{\ss} \ref{5.5} induziert eine abbildungsdefinierte quasi-analytische Zerlegung eine abbildungsdefinierte holomorphe Bl\"atterung. Im Fall von fasertreuen lokalen Beschreibungen gilt die Umkehrung.

\Beg{Satz}\label{5.15}
   Die $p$-dimensionale holomorphe Bl\"atterung $\5F$ sei abbildungsdefiniert mit fasertreuen 
   lokalen Beschreibungen. Dann gilt:
   \BegEN
      \item $\5F$ ist eine starke Bl\"atterung.\label{5.15.1}
      \item \label{5.15.2}$\DF$ ist abbildungsdefiniert mit den fasertreuen lokalen Beschreibungen von $\5F$ 
      als lokale Beschreibungen.
   \EndEN
\End

{\bf Beweis:}
Sei $U \subset X$ ein Gebiet und $\Mapping f U W$ eine fasertreue lokale Beschreibung von $\5F$. 
Es gen\"ugt zu zeigen,
da{\ss} die Elemente von $\5D_{f}$ starke lokale $\5F$-Bl\"atter sind: daraus folgt dann sofort \ref{5.15.1} und
\ref{5.15.2}.

Sei dazu $A \in \5D_{f}$. Wir betrachten einen Punkt $x \in A$. Dann gilt f\"ur alle Punkte $y$ einer $X$-offenen Umgebung $U' \subset U$ von $x$, da{\ss} $\dim_{y} f^{-1}\big(f(y)\big) \leq \dim_{x}A$ (vgl.~\cite{KK}).
Also ist $\dim_{x}A \geq p$. Sei wieder $x \in A$ und $g \in (\5O_{\5F})_{x}$. Dann ist $g|_{A}$ konstant wegen
\ref{5.10}. Also ist $A$ eine $\5F$-Integralvariet\"at mit $\dimMin A \geq p$.

Nun sei $B \subset U$ eine lokal-analytische Menge und eine $\5F$-Integralvariet\"at. Sei $x \in A \cap B$. Weil wir uns $W$ als lokal-analytische Teilmenge von $\C^N$ vorstellen d\"urfen, ist $f|_{B_{x}}$ konstant, also 
$B_{x} \subset A_{x}$. \Qed

Die abbildungsdefinierten quasi-analytischen Zerlegungen mit fasertreuen lokalen Beschreibungen entsprechen also genau den abbildungsdefinierten holomorphen Bl\"atterungen mit fasertreuen lokalen Beschreibungen.

\Beg{Satz}\label{5.16}
   Die quasi-analytische Zerlegung $\5D$ von $X$ sei reindimensional. Ist $\5D$ abbildungsdefiniert, so ist
   $\5D$ abbildungsdefiniert mit fasertreuen lokalen Beschreibungen, wird also insbesondere durch eine
   starke holomorphe Bl\"atterung induziert.
\End

Zum {\bf Beweis} siehe \ref{5.11}.

\Beg{Beispiel}
   Die quasi-analytische Zerlegung $\5D = \5D_{f}$ aus Beispiel \ref{5.14} kann nicht durch fasertreue lokale 
   Beschreibungen definiert werden.
\End  
   Nehmen wir an, da{\ss} dieses doch m\"oglich sei, so betrachten wir die 
   zugeh\"orige holomorphe Bl\"atterung $\5F$. Weil $\5F$ und $\5F_{f}$ auf $X \sm \Sing f$ \"ubereinstimmen, 
   mu{\ss} $\5F = \5F_{f}$ sein. Es ist $g \in \5O_{\5F}(X)$, aber $g$ ist auf dem Blatt $\{0\} \times \C^2$ von
   $\5D$ nicht konstant; Widerspruch!

\Beg{Beispiel}
   Wir fahren fort mit Beispiel \ref{5.14} und betrachten zum Vergleich die folgende Abbildung
   $\Mapping {\Dot f}XW$, $(z,w,u) \mapsto (z,w \cdot u)$. Sie hat die    Funktionalmatrix 
   $$D\Dot f = \pmatrix{1 & 0 & 0 \cr 0 & u & w},$$
   also ist $\Sing \Dot f = \C \times \{(0,0)\}$.
   F\"ur $(s,t) \in \C^2$ gilt:
   $$\Dot{f^{-1}}(s,t) = \{s\} \times \big\{(w,u) \in \C^2 : w \cdot u = t\big\}.$$
   Wir stellen fest:
   
   $\5D_{f}$ besteht aus den Bl\"attern $\{0\} \times \C^2$\  und \ 
   $A_{s,c} := \{s\} \times \big\{(w,u) \in \C^2 : w \cdot u=c\big\},   s \in \C^*, c \in \C$

   $\5D_{\Dot f}$ besteht aus den Bl\"attern
   $A_{s,c},   s \in \C, c \in \C$
   
   Au{\ss}erdem ist $\5D_{\Dot f}$ reindimensional, also ist $\Dot f$ fasertreu und
   $\5D_{\Dot f}|_{X \sm \Sing f} = \5D_{f}|_{X \sm \Sing f}$. Deshalb ist $\5F_{f} = \5F_{\Dot f}
    = :\5F$, insbesondere $\DF = \5D_{\Dot f}$. Es ist $\5D_{f} \neq \5D_{\Dot f}$; $\DF$
   kann nicht direkt aus $f$ berechnet werden.
\End

Wir schlie{\ss}en mit einem fundamentalen Ergebnis von {\sc Malgrange} (vgl.~Satz von Malgrange in \cite{Reiffen1})

\Beg{Satz}\label{5.19}
   Die Garbe $\Omega_{\5F}$ der $p$-dimensionalen holomorphen Bl\"atterung $\5F$ sei lokal frei und
   es gelte $\codim \Sing \5F \geq 3$. Dann ist $\5F$ abbildungsdefiniert mit fasertreuen lokalen Beschreibungen.
\End

\subsection{Koh\"arente, vollst\"andige und perfekte quasi-analytische Zerlegungen}\label{KZ}

Weil quasi-analytische Zerlegungen sehr chaotisch sein k\"onnen, ist es gewi{\ss} angebracht, an die Zerlegungen Forderungen zu stellen, die aus komplex-analytischer Sicht sinnvoll sind. Am Ende dieses Paragraphen findet man Beispiele, die unsere Begriffsbildungen motivieren und die Grenzen unserer Ergebnisse zeigen.

\Beg{Satz}\label{6.1}
   Sei $\5D$ eine rein $p$-dimensionale quasi-analytische Zerlegung von $X$. Ist $\Theta^{\5D}$ eine
   koh\"arente Untergarbe von $\Theta$, so sind folgende Aussagen \"aquivalent:
   \BegEN
      \item $\Sing \5D\neq X$,
      \item $\Rang \Theta^{\5D}  = p$.
   \EndEN
   Wenn diese Bedingungen erf\"ullt sind, dann ist $\Sing \5D = \Sing \Theta^{\5D}$.
\End
Zum {\bf Beweis} vgl.~\cite{Knoche}, Satz 7.

Wir verallgemeinern die in \ref{6.1} beschriebene Situation:

\Beg{Definition}\label{6.2}
   Die quasi-analytische Zerlegung $\5D$ von $X$ hei{\ss}t \DEFI{koh\"arent}\index{Zerlegung!koh\"arente}, wenn
   gilt:
   \BegEN
      \item $\Theta^{\5D}$ ist koh\"arent
      \item $\5D|_{X \sm \Sing \Theta^{\5D}}$ ist rein $p$-dimensional, wobei $p = \Rang \Theta^{\5D}$ sei.
   \EndEN
   Dann sei $p$ der \DEFI{Rang} von $\5D$; wir schreiben $\Rang \5D = p$.
\End

\Beg{Satz}\label{6.3}
   Die quasi-analytische Zerlegung $\5D$ sei koh\"arent. Dann ist $\Sing \5D = \Sing \Theta^{\5D}$.
\End

{\bf Beweis:} Wegen \ref{5.8} ist $\5D|_{X \sm \Sing \Theta^{\5D}}$ regul\"ar, also gilt 
$\Sing \5D \subset \Sing\Theta^{\5D}$. Sei $x \in X \sm \Sing \5D$. Dann ist $\5D$ in einer $X$-Umgebung von $x$ reindimensional, also rein $p$-dimensional, $p = \Rang \Theta^{\5D}$. Also mu{\ss} 
$x \in X \sm \Sing \Theta^{\5D}$ und damit $\Sing \Theta^{\5D} \subset \Sing \5D$ gelten. \Qed

Insbesondere gilt:

\Beg{}\label{6.4}
   Ist $\5D$ eine koh\"arente quasi-analytische Zerlegung vom Rang $p$, so ist $\Sing \5D$ eine niederdimensionale
   analytische Teilmenge von $X$ und $\5D|_{X \sm \Sing \5D}$ eine regul\"are $p$-dimensionale quasi-analytische 
   Zerlegung von $X \sm \Sing \5D$.
\End

\Beg{Definition}\label{6.5}
   Sei $\5D$ eine koh\"arente quasi-analytische Zerlegung vom Rang $p$. Die durch die Komplettierung 
   $\widetilde{\Theta^{\5D}}$ von $\Theta^{\5D}$ definierte $p$-dimensionale holomorphe Bl\"atterung $\5F^{\5D}$
   auf $X$ hei{\ss}e die \DEF{durch $\5D$ definierte Bl\"atterung}\index{Bl\"atterung!durch $\5D$ definiert}.
\End

In der Situation von \ref{6.5} ist $\5D|_{X \sm \Sing \5D}$ die zur regul\"aren Bl\"atterung 
$\5F_{\5D}|_{X \sm \Sing \5D}$ geh\"orige Zerlegung.

\Beg{}\label{6.6}
   Sei $\5D$ eine koh\"arente quasi-analytische Zerlegung. Dann gilt $\Theta^{\5D} \subset \Theta_{\5F^{\5D}}$.
\End

{\bf Beweis:} Es ist $\Theta_{\5F^{\5D}} = \widetilde{\Theta^{\5D}}$.

\Beg{Definition und Bemerkung}\label{6.7}
   Die koh\"arente quasi-analytische Zerlegung $\5D$ von $X$ hei{\ss}t \DEFI{voll\-st\"a\-ndig}, wenn
   $\Theta^{\5D}$ vollst\"andig ist, d.h.~$\Theta^{\5D} = \Theta_{\5F^{\5D}}$. In diesem Fall ist
   \BegEN
      \item $\Sing \5D = \Sing \5F^{\5D}$,\label{6.7.1}
      \item $\codim \Sing \5D \geq 2$.\label{6.7.2}
   \EndEN
\End

Wegen \ref{6.7.2} vgl.~\ref{6.7.1}.

\Beg{Satz}\label{6.8}
   Sei $\5D$ eine koh\"arente quasi-analytische Zerlegung von $X$, $\5F := \5F^{\5D}$. Dann sind folgende
   Aussagen \"aquivalent:
   \BegEN
      \item $\5D$ ist vollst\"andig.\label{6.8.1}
      \item $(X,\DF')$ ist eine quasi-analytische Teilmenge von $(X,\5D)$.\label{6.8.2}
   \EndEN
\End
\ref{6.8.2} bedeutet, da{\ss} $\DF'$ auf jedem $B \in \5D$ eine quasi-analytische Zerlegung definiert.

{\bf Beweis:}
Ò\ref{6.8.2} $\IfThen$ \ref{6.8.1}Ó: Es ist $\Theta^{\5D} \subset \Theta_{\5F}$. Sei $U \Opss X$ ein Gebiet und 
$\theta \in \Theta_{\5F}(U)$. Wir betrachten ein Blatt $B$ von $\5D|_{U}$. Sei $x \in B$ und $A$ das Blatt von 
$\DF'|_{U}$ mit $x \in A$. Wegen \ref{6.8.2} ist $A$ eine quasi-analytische Teilmenge von $B$. Dann folgt
$\theta|_{x} \in \TT(A,x) \subset \TT(B,x)$. Also ist $\theta \in \Theta^{\5D}(U)$ und 
$\Theta_{\5F} \subset \Theta^{\5D}$.

Ò\ref{6.8.1} $\IfThen$ \ref{6.8.2}Ó:  Sei $\0x \in X$ und $\0x \in A \cap B$, $A \in \DF'$, $B \in \5D$. Wir gehen von der im Beweis von \ref{4.9} f\"ur $\5D' = \DF'$ dargestellten Situation aus. Dabei d\"urfen wir annehmen, da{\ss} ein
$\5D$-Pl\"attchen $\0B$ von $B$ existiert mit
\BegIT
   \item $\0x \in \0B$,
   \item $\0B$ ist eine analytische Teilmenge von $\0X$.
\EndIT

Sei $\zeta \in \0A \subset \C^p_{z}$, $\zeta \neq 0$ und $\theta := \sum_{\nu=1}^p \zeta_{\nu}\theta^{(\nu)}$.
Dann gilt:
\BegIT
   \item $\theta \in \Theta_{\5F}(\0X) = \Theta^{\5D}(\0X)$,
   \item $\Mapping \gamma {\Clint 01}{\0A}$, $\gamma(t) := \zeta t$, ist eine Integralkurve zu $\theta$ mit
   $\gamma(0)=  0 = \0x$, $\gamma(1) = \zeta$. 
\EndIT
Weil $\theta \| \0B$, folgt Im $\gamma \subset \0B$, 
   $\zeta = \gamma(1) \in \0B$. Also ist $\0A \subset \0B$. \Qed

\Beg{Definition}
   Die quasi-analytische Zerlegung $\5D$ von $X$ hei{\ss}e \DEFI{fast-regul\"ar}, wenn $\5D$ koh\"arent ist und keine
   (bzgl.~der $\5D$-Struktur) irreduzible Komponente eines Blattes in $\Sing \5D$ enthalten ist.
\End

\Beg{Satz}\label{6.10}
   Jede fast-regul\"are quasi-analytische Zerlegung $\5D$ von $X$ ist reindimensional und voll\-st\"andig. 
\End

{\bf Beweis:} Die Aussage \"uber die Reindimensionalit\"at ist klar. Sei $U \Opss X$, $\theta \in \Theta_{\5F^{\5D}}(U)$, $A \in \5D$,
$A \cap U \neq \emptyset$. Dann ist $\theta$ auf $U \sm \Sing \5D$ parallel zu $A$. Wegen \ref{1.14} gilt dies auf ganz $U$. \Qed

\Beg{Satz}\label{6.11}
   F\"ur die quasi-analytische Zerlegung $\5D$ von $X$ gelte:
   \BegEN
      \item $\Sing \5D$ ist eine analytische Teilmenge von $X$ mit $\codim \Sing \5D \geq 2$.
      \item Keine (bezgl.~der $\5D$-Struktur) irreduzible Komponente eines Blattes ist in $\Sing \5D$ enthalten.    
   \EndEN
   Dann ist  $\5D$ fast-regul\"ar, insbesondere vollst\"andig.
\End

{\bf Beweis:} Sei $S := \Sing \5D$. $\Theta^{\5D}|_{X \sm S}$ ist eine regul\"are Untergarbe von 
$\Theta|_{X \sm S}$. Aufgrund eines Satzes von {\sc Siu-Trautmann} (vgl.~\cite[Prop. 1.21]{Reiffen1})
gibt es eine koh\"arente Untergarbe $\Theta'$ von $\Theta$ mit $\Theta'|_{X \sm S} = \Theta^{\5D}|_{X \sm S}$.  
Weil $\Theta^{\5D}|_{X \sm S}$ regul\"ar ist, d\"urfen wir $\Theta'$ sofort als vollst\"andig voraussetzen. Dann ist 
$\Theta^{\5D} \subset \Theta'$. Mit der Argumentation im Beweis von \ref{6.10} folgt $\Theta^{\5D} = \Theta'$.
\Qed

Aus \ref{6.11} folgt:

\Beg{Korollar}\label{6.12}
   Die quasi-analytische Zerlegung $\5D$ von $X$ sei rein 1-codimensional und $\Sing \5D$ sei eine analytische 
   Teilmenge von $X$ der Codimension $\geq 2$. Dann ist $\5D$ fast-regul\"ar, insbesondere vollst\"andig.
\End

\Beg{Satz}\label{6.13}
   Sei $\5F$ eine starke holomorphe Bl\"atterung. Dann gilt:
   \BegEN
      \item $\5D := \DF$ ist vollst\"andig,\label{6.13.1}
      \item $\5F^{\5D} = \5F$.\label{6.13.2}
   \EndEN
\End

{\bf Beweis} ad \ref{6.13.1}: 
Wegen \ref{4.16} ist $(X,  \DF')$ eine quasi-analytische Teilmenge von $(X,\DF)$. Weil $\Theta_{\5F}$ vollst{\"a}ndig ist, folgt: $\Theta^{\5D}\subset\Theta_{\5F}$. Mit dem Schlu{\ss} \ref{6.13.2} $\IfThen$ \ref{6.13.1} im Beweis von \ref{6.8} folgt: 
$\Theta^{\5D} =\Theta_{\5F}$. 
\ref{6.13.2} ist wegen \ref{6.13.1} klar. \Qed

Ist $\5F$ eine starke holomorphe Bl\"atterung, so ist $\5D_{\5F}$ insbesondere koh\"arent. Wir wissen nicht, ob dies f\"ur beliebige holomorphe Bl\"atterungen mit Bl\"attern \"uberall gilt. 

\Beg{}\label{6.14}
   Sei $\5F$ eine holomorphe Bl\"atterung mit Bl\"attern \"uberall und $\5D := \DF$ sei koh\"arent. Dann gilt
   $\5F^{\5D} = \5F$.
\End

{\bf Beweis:} Es gilt mit bekannten Argumenten $\widetilde{\Theta^{\5D}} = \Theta_{\5F}$.

\Beg{Definition}\label{6.15}
   Die quasi-analytische Zerlegung $\5D$ von $X$ hei{\ss}t \DEFI{perfekt}, wenn gilt:
   \BegEN
      \item $\5D$ ist koh\"arent.
      \item Es gibt eine holomorphe Bl\"atterung $\5F$ mit Bl\"attern \"uberall derart, da{\ss} $\5D = \DF$ ist.
   \EndEN
\End

Die perfekten quasi-analytischen Zerlegungen sind also genau die koh\"arenten Bl\"atterr\"aume der Bl\"atterun\-gen mit Bl\"attern \"uberall.

Aus \ref{6.14} folgt:

\Beg{}\label{6.16}
   Die quasi-analytische Zerlegung $\5D$ sei perfekt und sei $\5F$ wie in \ref{6.15}. Dann ist $\5F = \5F^{\5D}$.
   Insbesondere ist $\5F$ eindeutig bestimmt.
\End

Eine alternative Definition der Perfektheit lautet:

\Beg{}\label{6.17}
   Die koh\"arente quasi-analytische Zerlegung $\5D$ ist genau dann perfekt, wenn gilt:
   $\5F:=\5F^{\5D}$ hat Bl\"atter \"uberall, $\DF$ ist koh\"arent und $\DF = \5D$.
\End

Aus \ref{5.5}, \ref{5.15} und \ref{6.13} folgt:

\Beg{Satz}\label{6.18}
   Die quasi-analytische Zerlegung $\5D$ sei abbildungsdefiniert mit fasertreuen lokalen Beschreibungen. 
   Dann ist $\5D$ vollst\"andig und perfekt.
\End

\Beg{Satz}\label{6.19}
   F\"ur die koh\"arente quasi-analytische Zerlegung $\5D$ gelte:
   \BegEN
      \item \label{6.19.1}$\5D$ ist rein $p$-dimensional und $\dim \Sing \5D < p$,
      \item \label{6.19.2}$\5D$ ist lokal eigentlich.
   \EndEN
   Dann ist $\5D$ vollst\"andig und perfekt.
\End

{\bf Beweis:} Wegen \ref{6.19.1} ist $\5D$ fast-regul\"ar und wegen \ref{6.10} dann vollst\"andig. 
Wegen \ref{6.7} ist $\Sing \5D = \Sing \5F_{\5D}$. Sei $\5F := \5F^{\5D}$, $S := \Sing \5F$.
Wir d\"urfen annehmen, da{\ss} $\5D$ eine analytische Zerlegung von $X$ ist. Auf $X^* := X \sm S$ stimmt $\5D$
mit dem Bl\"atterraum $\5D_{\5F^*}$ der regul\"aren Bl\"atterung $\5F^* := \5F|_{X^*}$ \"uberein. Ist $A'$ eine irreduzible Komponente eines Blattes $A$ von $\5D$, so ist $A' \sm S$ ein Blatt von $\5F^*$. Ist umgekehrt $A^*$ ein Blatt von $\5F^*$, so ist die $X$-abgeschlossenen H\"ulle $A'$ von $A^*$ eine irreduzible Komponente eines Blattes $A$ von $\5D$.

Sei $A$ ein Blatt von $\5D$. Dann ist $A$ eine Integralvariet\"at von $\5F$ und $\dimMin A \geq p$. Sei $B \subset X$ eine irreduzible lokal-analytische Menge, Integralvariet\"at von $\5F$, $\dimMin B \geq p$ und $x \in A \cap B$. Wir wollen zeigen: $B_{x} \subset A_{x}$. Indem wir uns eventuell auf eine Umgebung von $x$ zur\"uckziehen, d\"urfen wir annehmen, da{\ss} $B$ eine analytische Teilmenge von $X$ ist. Dann ist $B \sm S$ ein Blatt von $\5F^*$, also ist $B$ eine irreduzible Komponente eines Blattes von $\5D$. Es folgt: $B \subset A$. Also ist $\5D = \DF$ perfekt. \Qed

\ref{6.19} kann auch folgenderma{\ss}en formuliert werden:

\BegNN{\ref{6.19}$'$ Satz}\label{6.19'}
   F\"ur die quasi-analytische Zerlegung $\5D$ gelte:
   \BegEN
      \item \label{6.19'.1}
         $\Sing \5D$ ist eine analytische Teilmenge von $X$, $\5D$ ist rein $p$-dimensional und $\dim \Sing \5D < p$,
      \item \label{6.19'.2}$\5D$ ist lokal eigentlich.
   \EndEN
   Dann ist $\5D$ vollst\"andig und perfekt.
\EndNN

{\bf Beweis:} Wir d\"urfen $p \leq n-1$ und dann $\codim \Sing \5D \geq 2$ voraussetzen. Wegen \ref{6.11} ist 
$\5D$ dann insbesondere koh\"arent, so da{\ss} die Voraussetzungen von \ref{6.19} erf\"ullt sind.

F\"ur die Dimension $p = n-1$ gibt es eine \"uber \ref{6.19} bzw.~\ref{6.19}$'$ hinausgehende Aussage:

\Beg{Satz}\label{6.20}
   F\"ur die koh\"arente quasi-analytische Zerlegung $\5D$ gelte:
   \BegEN
      \item $\5D$ ist rein 1-codimensional,\label{6.20.1}
      \item $\5D$ ist lokal eigentlich.
   \EndEN
   Dann ist $\5D$ vollst\"andig und perfekt.
\End

{\bf Beweis:} Wir werden zeigen, da{\ss} in diesem Fall automatisch $\codim \Sing \5D \geq 2$ ist; dann folgt
\ref{6.20} aus \ref{6.19}.

Wir d\"urfen voraussetzen, da{\ss} $\5D$ analytisch ist. Sei $\5F := \5F^{\5D}$, $S := \Sing \5D$, $X^* := X \sm S$ und $S' := \Sing \5F$. Es ist $\codim S' \geq 2$.

Wir gehen indirekt vor und nehmen an, da{\ss} ein $\0x \in S$ mit $\dim_{\0x}S = n-1$ existiert. 
Weil $\codim S'\geq2$ und $\codim \Sing S \geq 2$ ist, d\"urfen wir annehmen, da{\ss} 
$\0x \in S \sm (S' \cup \Sing S)$ ist. Deshalb k\"onnen wir von der folgenden Situation ausgehen:
\BegIT
   \item 
      $X$ ist ein Gebiet im $\C^n$ von der Form $X  = \1D \times \2D$, wobei $\1D$ ein Gebiet im $\C^{n-1}_{z}$
      und $\2D$ ein Gebiet in $\C_{w}$ ist, es ist $\0x = 0$.
   \item
      $\5F$ ist eine regul\"are Bl\"atterung mit dem Bl\"atteraum $\5D^* = \big\{ D_{1} \times \{w\} : w \in \2D\big\}$
   \item
      $S$ ist eine zusammenh\"angende $1$-codimensionale Untermannigfaltigkeit von $X$.
\EndIT
Wir produzieren den gew\"unschten Widerspruch, indem wir zeigen: $\5D^* = \5D$. Dazu unterscheiden wir zwei alternative F\"alle:

{\bf 1.~Fall:}
Die Bl\"atter von $\5D^*$ schneiden $S$ alle niederdimensional. 
\\
Weil $\5D^*|_{X^*} = \5D|_{X^*}$ ist, mu{\ss} dann $\5D^* = \5D$ sein. 

{\bf 2.~Fall:}
Es gibt ein Blatt von $\5D^*$, welches $S$ volldimensional schneidet. 
\\
Dann mu{\ss} dieses Blatt gleich $S$ sein, d.h.~es ist $S = \2D \times \{\0c\}$ f\"ur ein $\0c \in \C_{w}$.
Weil $0 \in S$ gilt, ist dann $\0c = 0$. Sei $A$ das Blatt von $\5D$ mit $0 \in A$. Weil $A$ nicht in $X^*$ eindringen kann, mu{\ss} $A = S$ sein. Also ist auch in diesem Fall $\5D^* = \5D$. \Qed

\Beg{Satz}\label{6.21}
   Die quasi-analytische Zerlegung $\5D$ sei  rein 1-codimensional und lokal eigentlich. Dann sind folgende
   Aussagen \"aquivalent:
   \BegEN
      \item\label{6.21.1}
         $\5D$ ist koh\"arent.
      \item\label{6.21.2}
         $\5D$ ist vollst\"andig.
      \item\label{6.21.3}
         $\5D$ ist perfekt.
      \item\label{6.21.4}
         $\5D$ ist abbildungsdefiniert.
    \EndEN
    Wegen \ref{5.11} ist  $\5D$ dann abbildungsdefiniert mit fasertreuen lokalen Beschreibungen. 
\End

{\bf Beweis:} Trivialerweise gilt \ref{6.21.2} $\IfThen$ \ref{6.21.1} und \ref{6.21.3} $\IfThen$ \ref{6.21.1}; wegen
\ref{6.20} gilt \ref{6.21.1} $\IfThen$ \ref{6.21.2}, \ref{6.21.1} $\IfThen$ \ref{6.21.3}. Also sind
\ref{6.21.1},  \ref{6.21.2} und \ref{6.21.3} \"aquivalent.

Wegen \ref{6.18} gilt \ref{6.21.4} $\IfThen$ \ref{6.21.3}. Es bleibt \ref{6.21.3} $\IfThen$ \ref{6.21.4} zu zeigen. Sei dazu $\5F := \5F^{\5D}$. Wegen \ref{6.21.2} ist $S := \Sing \5D = \Sing \5F$. Dann sind die Voraussetzungen des Satzes von {\sc Mattei-Moussu} (vgl.~3.29 in \cite{Reiffen1}) f\"ur $\5F$ erf\"ullt. Deshalb ist $\5F$ abbildungsdefiniert mit holomorphen Funktionen (die sind stets offene Abbildungen) als lokalen Beschreibungen. Insbesondere ist 
$\5F$ eine starke Bl\"atterung. 

Sei $\Mapping f U \C$ eine lokale Beschreibung von $\5F$. Dann ist $\5D_{f} = \DF|_{U}$ (vgl. \ref{5.15}). Wir d\"urfen $\5D|_{U}$ als analytisch voraussetzen. Auf $U \sm S$ ist $\5D = \DF$. Weil $\codim S \geq 2$ ist, mu{\ss} $f$ auf den Bl\"attern von $\5D|_{U}$ konstant sein. Dann ist aber $\5D|_{U} = \5D_{f}$. \Qed

F\"ur perfekte quasi-analytische Zerlegungen gilt ein Identit\"atssatz:

\Beg{Satz}\label{6.22}
   Seien $\5D, \Dot{\5D}$ perfekte quasi-analytische Zerlegungen von $X$. Gibt es ein nichtleeres Gebiet 
   $U \Opss X$ mit $\5D|_{U} = \Dot{\5D}|_{U}$, so ist $\5D = \Dot{\5D}$.
\End

Klar! Es ist $\5F^{\5D} = \5F^{\Dot{\5D}}$.

Es folgen einige Beispiele:

\Beg{Beispiel} (vgl.~\cite[Beispiel 15]{Knoche})
  Es seien $Y$ und $Z$ komplexe zusammenh\"angende Mannigfaltigkeiten 
mit
abz\"ahl\-barer Topologie von der gleichen Dimension $m$, es sei 
$X:=Y\times Z$. Ferner
sei $U \Opss Y$ nicht leer, offen und zusammenh\"angend, es sei $A:=
Y\backslash U$. Dann wird durch
$$\5D(y,z):=\left\{
\begin{array}{ll}
U\times \{z\} &\hbox{ wenn } y\in U,\\
\{y\}\times Z &\hbox{ wenn } y\notin 
U\end{array}\right. \hbox{  f\"ur }
(y,z)\in X
$$
eine rein $m$-dimensionale lokal-analytische Zerlegung $\5D$ 
auf $X$
definiert. Es ist $\Sing \5D=\partial U\times Z$; die Garbe 
$\Theta^\5D$ ist
genau dann koh\"arent, wenn $A$ analytisch in $Y$ ist (dann ist 
$A=\partial U$):
wenn $\Theta^{\5D}$ koh{\"a}rent ist, dann ist $\Sing\5D = \partial U \times Z$ sowie auch $A = \partial U$ analytisch; wenn umgekehrt $A$ analytisch ist, dann
ist mit der  Idealgarbe von $A$ auch $\Theta^{\5D}$ koh{\"a}rent.

Wenn $\Theta^\5D$ koh\"arent und $A\neq \emptyset$ ist, dann ist 
$\Theta^\5D$
sicher nicht vollst\"andig, denn die Komplettierung
$\5F:=\widetilde{\Theta^\5D}$ ist diejenige regul\"are holomorphe 
Bl\"atterung
von $X$, deren Bl\"atter genau die $Y\times \{z\},\, z\in Z$, sind 
(vgl.~Identit\"atssatz f\"ur singul\"are Bl\"atterungen). Folglich ist 
$\Theta^\5D$ genau
dann vollst\"andig, wenn $A=\emptyset$ ist. Wenn $A$ 
innere Punkte hat, ist
$\Theta^\5D$ nicht koh\"arent. Ist $A\neq\emptyset$, so ist $\5D$ nicht
eigentlich.

Ist z.B.~$Y=Z=\C$ und $U=\C^\ast$, dann ist $\5D$ rein 
1-codimensional,
$\Theta^\5D$ wird erzeugt vom Vektorfeld $y\frac{\partial}{\partial y}$, 
ist
also koh\"arent, aber nicht vollst\"andig: $\widetilde{\Theta^\5D}$
wird erzeugt
von $\frac{\partial}{\partial y}$; es ist
$\dim\Sing\5D=1=\dim\5D$. 
\End

\Beg{Beispiel} (vgl.~\cite[Beispiel 17]{Knoche})
 Es sei $X:= \C^3$, dann wird durch
$$\5D(z):=\left\{
\begin{array}{ll}
\C\times \{z_2\}\times \{z_3\} &\hbox{ wenn $z_3$  
rationalen Real-
und Imagin\"arteil hat}\\
\{z_1\}\times\C\times \{z_3\}  &\hbox{ sonst } \end{array}\right. 
$$
eine rein 1-dimensionale analytische Zerlegung von $X$ definiert; 
die Garbe
$\Theta^\5D=0\subset \Theta$ definiert eine regul\"are Bl\"atterung 
mit dem
einzigen Blatt $X$, es ist $X=\Sing \5D\neq\Sing\Theta^\5D=\emptyset$.
$\5D$ ist eigentlich. 
\End

\Beg{Beispiel} (vgl.~\cite[Beispiel 16]{Knoche})
   
F\"ur die durch
$$\5D(y,z):=\left\{
\begin{array}{ll}
\{(z_1,z_2)\} \times \C &\hbox{  wenn } z_1\in \C^\ast,\\
\{0\}\times \C\times \{z_3\} &\hbox{  wenn } z_1=
0\end{array}\right.
$$
auf $\C^3$ definierte rein 1-dimensionale analytische Zerlegung wird
$\Theta^\5D$ erzeugt vom globalen Vektorfeld 
$z_1\frac{\partial}{\partial
z_3}$, ist also koh\"arent; die  zugeh{\"o}rige Bl\"atterung $\5F$ wird erzeugt vom globalen Vektorfeld 
$\frac{\partial}{\partial
z_3}$; die Bl\"atter von $\5F$ sind also genau die Mengen 
$\{(z_1,z_2)\}\times
\C$, wobei $(z_1,z_2)\in\C^2$. Es ist
$$
\Sing \5D=\{z\in\C^3:z_1=0\}=\Sing \Theta^\5D,\quad \Sing\5F =\emptyset.
$$
$\5D$ ist eigentlich. 
\End

\Beg{Beispiel}\label{6.26}
   Sei $X = \C^2_{z,w}$ und $\5D$ bestehe aus den Mengen
   $$A = \{0\} \times \C_{w}, \qquad A_{a} = \{(z,w) : w = az, z \neq 0\}, \ a \in \C.$$
   $\5D$ ist eine rein 1-codimensionale lokal-analytische Zerlegung von $X$. $\Theta^{\5D}$ wird erzeugt von
   $\theta := z\,\frac \partial {\partial z} + w\,\frac \partial{\partial w}$, folglich ist $\Theta^{\5D}$  koh\"arent. Weil $\Sing \5D = \{0\}$ 
   2-codimensional ist, ist $\Theta^{\5D}$ vollst\"andig. Die Bl\"atterung $\5F := \5F^{\5D}$ besitzt kein durch 0 verlaufendes Blatt. $\5D$ ist nicht lokal eigentlich.
\End

\subsection{Anhang}\label{Anh}

 \Beg{Definition}
    Ein topologischer Raum $X$ hei{\ss}t \DEFI{(abgeschlossen) zerlegbar}\index{zerlegbar},
    wenn $X$ eine Darstellung 
    $$X = \bigcup_{\nu \in N} A_{\nu}\eqno(*)$$
     besitzt, wobei $N$ eine Teilmenge von $\N$
    mit {\rm Anz} $N \geq 2$ ist und die $A_{\nu}$ nichtleere paarweise disjunkte abgeschlossene Teilmengen 
    von $X$ sind.
    Das System $\{A_{\nu}:\nu \in N\}$ hei{\ss}t dann eine \DEFI{(abgeschlossene) Zerlegung}\index{Zerlegung}
     von $X$. Wenn man die
    Menge $N$ endlich w\"ahlen kann, dann hei{\ss}t $X$ 
    \DEFI{endlich (abgeschlossen) zerlegbar}.
    
    Wenn $X$ hausdorffsch ist, dann nennen wir $X$ \DEFI{kompakt zerlegbar} 
    bzw.~\DEFI{endlich kompakt zerlegbar}, wenn alle $A_{\nu}$ kompakt gew\"ahlt werden k\"onnen. In dem Fall
    nennen wir die zugeh\"orige Zerlegung eine \DEFI{kompakte Zerlegung}.
 \End
 
 \BegIT
    \item
      Ein topologischer Raum  ist genau dann endlich zerlegbar, wenn er nicht zusammenh\"angend ist. 
   \item
      Ein kompakter Raum ist genau dann kompakt zerlegbar, wenn er zerlegbar ist.
\EndIT
      
Uns ist nicht bekannt, ob ein beliebiger zusammenh{\"a}ngender Raum nicht zerlegbar ist.

Allerdings k\"onnen wir zeigen:
 
 \Beg{Satz}\label{7.2}
    Ein wegzusammenh\"angender topologischer Raum ist nicht zerlegbar. 
 \End
 
 {\bf Beweis:} Wir \"uberlegen uns zun\"achst, da{\ss} es gen\"ugt, den Spezialfall $X=\Clint 01$ zu behandeln: angenommen, es gibt eine Zerlegung $\{A_{\nu}:\nu \in N\}$ von $X$. 
W\"ahle $\nu \ne \mu$ aus $N$, $a \in A_{\nu}$, $b \in A_{\mu}$ und einen Weg $\Mapping \ga {\Clint 01}X$ von $a$ nach $b$. Dann kann man aus den kompakten Mengen $\ga^{-1}(A_{\la}) \subset \Clint 01$, $\la \in N$, eine kompakte Zerlegung von $\Clint 01$ konstruieren.

Angenommen, es gibt eine kompakte Zerlegung $\{A_{\nu}:\nu \in \N\}$ von $\Clint 01$. Wir d\"urfen annehmen, dass $0 \in A_{0}$ und $1 \in A_{1}$. Es sei $a := \max A_{0}$ und $b := \min \big(A_{1} \cap \Opclint a1\big)$. 
Man \"uberlegt sich, da{\ss} man $a < b$ annehmen darf. Dann ist
$I := \Opint ab$ ein offenes Intervall und $\{A_{\nu} \cap I : \nu \geq 2\}$ eine kompakte Zerlegung von $I$.
Wegen Lemma \ref{7.3} d\"urfen wir annehmen, da{\ss} 
alle $K_{\nu} := A_{\nu}$ kompakte Intervalle in $I$ sind (dabei sind entartete Intervalle, die nur aus einem Punkt bestehen, zugelassen). Es sei $M_{\nu}$ das Innere von $K_{\nu}$; dann ist $M_{\nu}$ genau dann leer, 
wenn $K_{\nu}$ nur aus einem Punkt besteht. Die Menge  $M := \bigcup_{\nu \geq 2}M_{\nu}$ ist offen, also ist 
$D :=\Clint ab  \sm M$ kompakt.  Da $\partial K_{\nu}$, $\nu \geq 2$,  jeweils aus h\"ochstens zwei Punkten besteht, 
ist $D$ abz\"ahlbar. Es ist $\dist(K_{\nu},K_{\mu})>0$ f\"ur $\nu \not= \mu$, $\nu,\mu \geq 2$. Daraus folgert man: 
$D$ ist perfekt, d.h.~jeder Punkt von $D$ ist H\"aufungspunkt von $D$, und $D$ ist total unzusammenh\"angend. Deshalb ist $D$ hom\"oomorph zur Cantormenge (vgl.~Abschnitt 30 in \cite{Willard}), insbesondere \"uberabz\"ahlbar. Widerspruch! \Qed

\Beg{Lemma}\label{7.3}
   Es sei $(A_{\nu})_{\nu \in \N}$ eine Folge paarweise disjunkter kompakter Teilmengen eines offenen Intervalles
   $I \subset \R$. Dann
   gibt es eine Folge $(L_{\nu})_{\nu \in \N}$ paarweise disjunkter
   (m\"oglicherweise leerer) kompakter Teilmengen von $I$ derart, dass
   f\"ur alle $n \in \N$ gilt:
   \BegEN
      \item $L_{n}$ ist endliche Vereinigung von paarweise disjunkten kompakten Intervallen.
      \item $\bigcup_{\nu \leq n}A_{\nu} =:{\bf A}_{n}   \subset {\bf L}_{n}:=\bigcup_{\nu \leq n}L_{\nu}$
      \item $\partial {\bf L}_{n} \subset {\bf A}_{n}$
   \EndEN
\End

{\bf Beweis} durch Induktion \"uber $n$: Setze ${L}_{0}:=\Clint{\min A_{0}}{\max A_{0}}$. Wenn $L_{0},\ldots,L_{n}$ mit den obigen Eigenschaften schon definiert sind, dann konstruieren wir $L_{n+1}$ folgendermassen: die
offene  Menge $ I \sm {\bf L}_{n}$ ist endliche Vereinigung von paarweise disjunkten offenen Intervallen $J_{\nu}^{(n)}$, $1 \leq \nu \leq j_{n}$. Die Mengen
$A^{(n+1)}_{\nu}:= A_{n+1} \cap J_{\nu}^{(n)}$, $1 \leq \nu \leq j_{n}$, sind kompakt wegen (3), und
$$L_{n+1}:= \bigcup_{1 \leq \nu \leq j_{n},\ A^{(n+1)}_{\nu}\not=\emptyset}\Clint{\min A^{(n+1)}_{\nu}}{\max A^{(n+1)}_{\nu}}$$
ist endliche Vereinigung paarweise disjunkter kompakter Intervalle.
Da $A_{n+1} \subset {\bf L}_{n} \cup L_{n+1} = {\bf L}_{n+1}$, folgt ${\bf A}_{n+1} \subset {\bf L}_{n+1}$.
Nach Konstruktion ist
$\partial L_{n+1} \subset  A_{n+1}$, also auch $\partial {\bf L}_{n+1} \subset {\bf A}_{n+1}$. \Qed

Bei den nun folgenden ¬\"Uberlegungen sei $a \in \R$ irrational und 
$$\AA:= \{ak + l : k,l \in \Z\}.$$

 Es gilt:
 
\Beg{Satz}\label{7.4}
   $\AA$ ist dicht in $\R$.
\End

Der {\bf Beweis} ist eine einfache Aufgabe f\"ur einen Anf\"angerkurs in Analysis.

Der Torus $T := \C/\Z^2$ \big(mit $(x+iy)\sim(x'+iy') \iff (x-x', y-y')\in \Z^2$\big) ist in nat\"urlicher Weise eine kompakte Riemannsche Fl{\"a}che, die kanonische Projektion $\Mapping\pi{\C}T$ ist lokal biholomorph. 
Wir betrachten die Gerade $\Mapping \gamma \R T$, $\gamma(t):=t + iat$. Bei der folgenden reellen Betrachtung sehen wir $\C$ als $\R^2$ an.

\Beg{Satz}\label{7.5}
   Es sei 
   $\xi \in T$ und $\eps \in \R$ mit $\eps > 0$. Dann gibt es $(x,y) \in \R^2$ mit 
   $\pi(x,y)  = \xi$ und $|ax - y| < \eps$.
\End
{\bf Beweis:}
 Es sei $(u,v) \in \R^2$ mit $\pi(u,v) = \xi$. Wegen \ref{7.4} gibt es $k,l \in \Z$ mit $|(au-v)-(ak+l)|<\eps$.
Dann haben $x:=u-k$ und $y:=v+l$ die geforderten Eigenschaften.

\Beg{Satz}\label{7.6}
   Die Menge $\pi\big( \{t + i\,at: t \in \R\}\big)$ ist dicht in $T$.
\End

{\bf Beweis:} Zu gegebenem $\xi \in T$ und $n \in \N$ mit $n>0$ gibt es wegen \ref{7.5} Zahlen $t_{n}$
und $y_{n} \in \R$ derart, dass $\xi = \pi(t_{n},y_{n})$ und $|at_{n}-y_{n}|< 1/n$ f\"ur alle $n$. Dann ist 
$\xi = \lim\limits_{n \to \infty}\pi(t_{n} + i\,at_{n})$. \Qed

Der Beweis der folgenden Aussage ist  eine \"Ubungsaufgabe f\"ur einen Kurs in Funktionentheorie:

\Beg{Satz}\label{7.7}
  Die Funktionen $\Mapping \sin \C\C$ und $\Mapping \cos \C \C$ sind surjektiv.
\End

{\tt burchard.kaup@unifr.ch}\\
D\'epartement de math\'ematiques, Universit\'e de Fribourg, CH-1700 Fribourg, Switzerland

{\tt reiffen@mathematik.uni-osnabrueck.de}\\
Fachbereich Mathematik/Informatik, Universit\"at Osnabr\"uck, D-49069 Osnabr\"uck, Germany

 \end{document}